\documentclass[final]{autart}

\usepackage{natbib}

\usepackage{amsmath, amssymb, amsfonts}
\usepackage[amsmath,thmmarks,amsthm]{ntheorem}
\usepackage{enumerate}
\usepackage{graphicx}
 
\usepackage{verbatim}
\usepackage{url}
\usepackage{kbordermatrix}

\usepackage{graphicx}
\usepackage{color}

\usepackage{amsmath}
\usepackage{multicol}

\newcommand{\trace}{\operatorname{trace}}

\newcommand*\diff{\mathop{}\!\mathrm{d}}

\newcommand\mycom[2]{\genfrac{}{}{0pt}{}{#1}{#2}}

\renewcommand{\kbldelim}{(}
\renewcommand{\kbrdelim}{)}

\newtheorem{Theorem}{Theorem}
\newtheorem{Assumption}{Assumption}
\newtheorem{Proposition}{Proposition}
\newtheorem{Lemma}{Lemma}
\newtheorem{Corollary}{Corollary}
\newtheorem{Definition}{Definition}
\newtheorem{Example}{Example}

\newtheorem{Remark}{Remark}

\newcommand {\R}{\mathbb R}

\newcommand {\C}{\mathbb C}
\newcommand {\M}{\mathbb M}

\newcommand {\I}{\mathbb I}

\newcommand{\be}{\begin{equation}}
\newcommand{\ee}{\end{equation}}
\newcommand{\sgn}{\operatorname{{\mathrm sgn}}}
\newcommand{\V}{\mathcal V}

\newcommand{\W}{\mathcal W}

\newcommand{\updt}[1]{{\color{blue}#1}}

\newcommand*\dif{\mathop{}\!\mathrm{d}}

\begin{document}

\begin{frontmatter}
\title{Revisiting Totally Positive Differential Systems: A~Tutorial and New Results}
\thanks[footnoteinfo]{
Research  supported in part 
by  research grants from    
   the Israel Science Foundation,  the 
 US-Israel Binational Science Foundation, 
 DARPA FA8650-18-1-7800, NSF 1817936, and AFOSR FA9550-14-1-0060.  A highly  abridged version of this paper
has been accepted for presentation at   the \emph{57th IEEE CDC}~\citep{tpds_cdc}.}
 \author[First]{Michael Margaliot}
\author[Second]{Eduardo D. Sontag}
\address[First]{M. Margaliot (Corresponding Author) is with the Department of  Elec. Eng.-Systems,  Tel-Aviv
 University, Tel-Aviv 69978, Israel. E-mail:  \texttt{michaelm@eng.tau.ac.il}}
\address[Second]{E. D. Sontag  is with the   Dept. of Elec. and Computer Eng. and the Dept. of
Bioengineering,  Northeastern University,
Boston, MA.
E-mail: \texttt{e.sontag@northeastern.edu}
}

\maketitle

\begin{abstract}
A matrix is called totally nonnegative~(TN)  [totally positive~(TP)] if all its minors
 are nonnegative [positive]. Multiplying a vector by a~TN matrix does not    increase  the number of
sign variations in the vector. In a largely forgotten paper,  Binyamin \cite{schwarz1970}
considered matrices whose \emph{exponentials}  are~TN or~TP. He also
  analyzed the evolution of the number of sign changes
in the vector solutions of the corresponding  linear system. 
His work, however, considered only linear systems. 

In a seemingly different line of research, 
  \cite{smillie}, \cite{periodic_tridi_smith}, and others analyzed 
 the stability of nonlinear tridiagonal cooperative systems
 by using  the number of sign  variations  in the derivative vector   as an integer-valued Lyapunov function.

 We   show   that these two research topics  are intimately
 related. This   allows   
to derive important  generalizations 
of    the results by \cite{smillie}  and \cite{periodic_tridi_smith} 
  while simplifying the proofs. These  generalizations are particularly relevant in the context of control systems. 
Also, the results by Smillie and Smith provide sufficient conditions for analyzing stability based
on the number of sign changes in the vector of derivatives and the connection to the work  of 
Schwarz allows to show  in what sense these results are also necessary. 
We   describe several   new and interesting  research directions
arising from this new connection.  

\end{abstract}
\end{frontmatter}

\emph{Keywords:}
 Totally nonnegative matrices,
 totally positive matrices, infinitesimal generators, 
cooperative  dynamical systems, 
monotone dynamical systems,   entrainment,  Floquet theory,
stability analysis, multiplicative and additive compound matrices. 
 
\section{Introduction}

The trajectories of monotone dynamical systems preserve a partial ordering, induced by a proper cone,
between their initial conditions. 
Hirsch's quasi-convergence theorem~\citep{hlsmith} shows that this property 
 has far reaching implications for the
 asymptotic behavior of their trajectories. 
An important special   case are cooperative systems arising when the cone that induces the partial ordering is the positive orthant. 
 In an interesting paper, \cite{smillie} considered the time-invariant, nonlinear, strongly cooperative, and  {tridiagonal} system:
\be\label{eq:fsys}
					\dot x(t)=f(x(t)),
\ee
where~$x(t) \in \R^n$ and~$f:  \R^n \to \R^n$,
and has shown that \emph{every}  trajectory either leaves any compact set or
converges to an equilibrium point. This result has found many applications as well as several interesting generalizations 
(see, 
e.g.~\cite{RFM_stability,chua_roska_1990,Donnell2009120,periodic_tridi_smith,fgwang2013}).\footnote{We note in passing 
that~\cite{Fiedler1999} have proved  a similar  result  using a very different technique.  
}

To explain Smillie's  proof,
let~$z:=\dot x$. Then~\eqref{eq:fsys} yields the variational equation
\be\label{eq:zdoteqjj}
				\dot z(t)= J(x(t))z(t),
	\ee
	where
$
				J(x):=\frac{\partial }{\partial x} f( x)
$
is the Jacobian of the vector field~$f$.
Smillie showed that  since~$J$ is tridiagonal with positive entries on the super- and sub-diagonal,
the number of sign variations in the vector~$z(t)$, denoted~$\sigma(z(t))$,
  is a non-increasing   function of~$t$. 
Recall  
	that for a vector~$y\in\R^n$ with no zero entries the number of sign variations
 in~$y$ is
\be\label{eq:sigmdfr}
\sigma(y):=\left|\{i \in \{1,\dots,n-1\} : y_i y_{i+1}<0\} \right | .
\ee
For example, for~$n=3$ consider  the vector~$z(\varepsilon):=\begin{bmatrix}  -1 & \varepsilon & 2 \end{bmatrix}' $.
 For \emph{any}~$\varepsilon  \in \R\setminus\{0\}$,  $\sigma(z(\varepsilon))$   is well-defined and equal to
one. More generally,
the function~$\sigma$    can be extended, via continuity, to the largest open set:
\begin{align*}
\V := & \{y\in\R^n:   y_1 \not =0,\; y_n \not=0,\; \text{and if }  y_i=0 
\\
&\text{ for some~$i \in \{2,\dots,n-1\} $ then } y_{i-1}y_{i+1}<0\}.
\end{align*}
  Note in particular that if~$y\in \V$ then~$y$ cannot have two adjacent zero coordinates. 
	
 To explain the basic idea underlying Smillie's proof, consider the case~$n=3$. Seeking a contradiction, assume that
	the sign pattern of~$z(t)$ near some time~$t=t_0$ is as follows:\\
\begin{center}
\begin{tabular}{c | c c c } 
& $t=t_0^-$ & $t=t_0$ & $t=t_0^+$ \\  \hline 
$z_1(t)$& $+$ & $+$ & $+$ \\  
$z_2(t)$ & $+$ & $0$ & $-$\\    
 $z_3(t)$& $+$ & $+$ & $+$   \\  
\end{tabular} \\
\end{center}
Note that in this case~$\sigma(t):=\sigma(z(t))$ \emph{increases} from~$\sigma(t_0^-)=0$ to~$\sigma(t_0^+)=2$.
However, using~\eqref{eq:zdoteqjj} and  the structure of~$J$  yields  
$
							\dot z(t_0)=\begin{bmatrix}  *& + & 0 \\ + & * & +\\ 0& + &* \end{bmatrix}
							\begin{bmatrix} +\\ 0 \\ +  \end{bmatrix},
$
where~$+ $ means a positive value, and~$*$ means some  value, and thus~$\dot z_2(t_0)> 0$, and the 
case described in the table above cannot take place. 
	Smillie's  analysis shows rigorously that when~$\sigma(t)$ changes it can only decrease.
	This is  based on direct
	analysis of the~ODEs and is non trivial 
	  due to the fact that
if an entry~$z_i(t)$ becomes zero at some time~$t=t_0$ (thus perhaps leading to a change in~$\sigma( t  )$ near~$t_0$)
one must consider the possibility that
	higher-order derivatives of~$z_i(t )$ are also zero at~$t=t_0$. 
	Smillie then used the behavior of~$\sigma(z(t))$ to deduce that
	for any point~$a$
	in the state-space of~\eqref{eq:fsys} the 
	  omega limit set~$\omega(a)$ 
		cannot include more than a single point, and thus every   trajectory either leaves any compact set or
	converges to an equilibrium point.

	\cite{periodic_tridi_smith}   extended Smillie's approach to the case
	of  a time-varying cooperative system  with a  tridiagonal Jacobian
	with positive entries
	on the super- and sub-diagonals for all time~$t$. He showed that 
	if the time-varying vector field is periodic with period~$T$
	then every solution of the nonlinear dynamical system 
	either leaves any compact set or
	converges to a periodic solution with period~$T$.
	\cite{fgwang2013}
	describe some generalizations of these ideas to time-recurrent systems. 
The main advantage of  this   analysis approach
 is that it allows to study 
  the global behavior of  
	the nonlinear system~\eqref{eq:fsys} 
	using the (time-varying) linear system~\eqref{eq:zdoteqjj}.

Here, we show that these results can be generalized, and their proofs simplified,
 by relating them  to a classical topic from linear algebra: the 
sign variation diminishing  property~(SVDP) of
totally nonnegative~(TN) matrices~\citep{total_book,pinkus}, and more precisely 
to the notion of \emph{totally positive differential systems}~(TPDSs) introduced by~\cite{schwarz1970}. 

  To explain this, we   recall   two more   definitions for   
		the number of sign variations in a vector~\citep{total_book}. 
 	For~$y\in\R^n$, let~$s^-(y)$ denote the number of sign variations
	in the vector~$y$ after deleting all its zero entries,
	and let~$s^+(y)$ denote the maximal possible number of sign variations
	in~$y$ after each zero entry is replaced by either~$+1$ or~$-1$. 
	Note that~$s^-(y) \leq s^+(y)$ for all~$y\in\R^n$. 
	For example, for~$y=\begin{bmatrix} 1& 0 &2 &-3 &0 &1.1 \end{bmatrix}'$, $s^-(y)=2$ and~$s^+(y)=4$. 
	 Let~$
			\W:=\{y\in\R^n:s^-(y)=s^+(y)\}.
$
Note   that  if~$y\in \W$ then~$y$ cannot have two adjacent zero coordinates.
An immediate yet important observation is that~$\W=\V$. Thus,
if~$y\in\W$ then~$s^-(y)=s^+(y)=\sigma(y)$.
	
	A classical result from the theory of~TN matrices~\citep{total_book}   states
	that if~$A\in \R^{n\times m}$ is totally positive~(TP) 
	and~$x \in \R^m \setminus \{0\}$ then~$s^+(Ax)\leq s^-(x)$,	
	whereas if~$A$ is~TN (and in particular if it is~TP) 
					then~$
					s^-(Ax)\leq s^{-}(x)$  for all~$x\in \R^m$. 
	To apply this~SVDP
	to the stability analysis of~\eqref{eq:fsys} note that if the \emph{transition matrix}
	corresponding to~$J(x(t))$ in~\eqref{eq:zdoteqjj} is~TP for all time~$t$ then we may expect   the number of sign variations in~$z(t)$ 
	to be a nonincreasing function of time. As already shown by \cite{smillie}, this implies that there exists a time~$s\geq 0$ such that~$x_1(t)$ and~$x_n(t)$ are monotone functions of time
	for all~$t\geq s$, and that every trajectory of~\eqref{eq:fsys}   either leaves any compact set
	or converges to an equilibrium point.

 \cite{schwarz1970} considered the linear matrix differential equation~$\dot Y(t)=A(t)Y(t)$, $Y(t_0)=I$,
with~$A(t)\in\R^{n\times n }$ a continuous matrix function of~$t$. Note that every  minor of~$I$ is nonnegative (as it is either zero or one). 
 Schwarz gave a formula for the induced 
dynamics of the minors of~$Y(t)$. He defined the system as a totally [nonnegative] positive dynamical system 
if for every~$t_0$ and every~$t>t_0$ the matrix~$Y(t)$ is TP~[TN].\footnote{We use here a slightly different terminology than that used  in~\cite{schwarz1970}.} His analysis is based on what is now known as
  the theory of cooperative dynamical systems: the system is a totally nonnegative dynamical system if the dynamics
 maps any set of nonnegative minors to a set of nonnegative  minors. 
However, the work of Schwarz  seems to have been largely  forgotten and its potential for the analysis of \emph{nonlinear}
dynamical systems has been overlooked. 
Indeed, according to Google Scholar Schwarz's paper  has been cited~$22$
times since its publication in~$1970$.

We show that    TPDSs  
can be immediately linked to recent results on the stability of nonlinear cooperative systems. 
We consider a   more general form than in~\cite{schwarz1970}, 
namely,~$\dot Y(t)=A(t)Y(t)$, with~$A(t)$ a 
measurable matrix  function of time rather than continuous.

 We then  show how this 
 can be used   
  to derive the  interesting  results of 
	\cite{smillie} and  \cite{periodic_tridi_smith} under milder technical
conditions and with simpler proofs.  These generalizations are particularly
 relevant to control systems.

	The next section reviews
	  relevant definitions and results  from the   theory of~TN matrices. Section~\ref{sec:main}
		reviews~TPDSs. Section~\ref{sec:app} shows how these results can be applied to analyze
		the stability of nonlinear time-varying
		tridiagonal cooperative  systems. 
		We believe that highlighting the deep connections between the work of
		Schwarz and more recent work on nonlinear tridiagonal cooperative systems
		opens the door for many new research directions. Some of these potential
		directions are described in Section~\ref{sec:nd}.
		The work of Schwarz has been largely forgotten, and~TN   matrices are not well-known in the systems and control community, so we 
  provide a self-contained tutorial   on   the tools  needed 
	for the stability analysis of nonlinear systems. 
	
 		We use standard notation. Vectors [matrices] are denoted by small [capital] letters. 
$\R^n$ is the set of vectors with~$n$ real coordinates.
For a (column) vector $x\in\R^n$, $x_i$ is the $i$th entry of $x$, and~$x'$ is the transpose of~$x$.
 Let~$\R^n_{++}:=\{v\in\R^n : v_i>0, i=1,\dots,n\}$, i.e. 
the set of all~$n$-dimensional vectors with positive entries. A square matrix~$B$ is called Metzler if every off-diagonal entry of~$B$ is nonnegative. 
The square identity matrix is denoted by~$I$, with dimension that should be clear from context.

	\section{TN and TP matrices}\label{sec:preli}
	
		We begin by  reviewing  known definitions and results from  the rich and beautiful
		theory of TN and TP 
		matrices that will be used later on. We consider only square and real matrices, 
		as this is the case that is relevant for our applications. For more information and proofs
		we refer to the two excellent monographs~\citep{total_book,pinkus}.
		Unfortunately, this field suffers from nonuniform terminology. We follow the more modern terminology in~\citep{total_book}.
		\begin{Definition}\label{def:tntp}
		A matrix~$A\in\R^{n\times n}$
		is called totally nonnegative [totally positive]
		if the determinant of \emph{every} square submatrix is nonnegative [positive]. 
		\end{Definition}
In particular, if~$A$ is~TN [TP] then every
		entry of~$A$ is nonnegative [positive]. 
		
\begin{Example}\label{exa:scale}
It is straightforward to verify that the   matrix
		$
		A=\begin{bmatrix} 0& 1&0 \\ 0 & 0& 1\\0&0&0 \end{bmatrix} 
		$ is~TN. Consider the 
		     permutation matrix~$P:=\begin{bmatrix}  0& 1& 0 \\ 0 & 0& 1 \\ 1& 0 &0 \end{bmatrix}$.
				Then~$ PAP^{-1}=\begin{bmatrix} 0& 1& 0\\ 0& 0& 0\\ 1& 0& 0 \end{bmatrix}$ 
	     is \emph{not}~TN because~$\det \left ( \begin{bmatrix}   0 & 1\\ 1& 0  \end{bmatrix} \right  ) <0 $. 
\end{Example}

		Determining  that an $n\times n$ matrix is TN [TP] by direct verification 
		of Def.~\ref{def:tntp}  
		requires checking    the signs of all its~$\sum_{k=1}^n  \binom{n}{k} ^2  $ minors. This is of course prohibitive. 
		Fortunately, the minors are not independent and thus there exist much more efficient methods for verifying
		that a matrix is~TN [TP]~\cite[Ch.~3]{total_book}. Also, 
		  some matrices with a special structure are known to be~TN.	We   review two such examples. The first is important 
				in  proving the SVDP of~TN matrices. 
				The second example, as we will see below, is closely related to Smillie's results. 
		\begin{Example}
		Let~$E_{i,j}$ denote the~$n\times n$ matrix with all entries zero, except for
		entry~$(i,j)$  that is one.
		For~$p \in \R$  and~$i\in\{2,\dots,n\}$, let
							\begin{align}\label{eq:eb}
							L_i(p):=I+p E_{i,i-1},\;\;
							U_i(p):=I+p E_{i-1,i}. 
							\end{align}
							Matrices in this form are called 
		\emph{elementary bidiagonal}~(EB) matrices. If the identity matrix~$I$
		in~\eqref{eq:eb} is replaced by a diagonal matrix~$D$
		 then the matrices are called
		\emph{generalized elementary bidiagonal}~(GEB).
    It is straightforward to see  that EB matrices 
		are TN when~$p\geq 0$, and that~GEB  matrices  are TN when~$p\geq 0$ and the diagonal matrix~$D$ is componentwise nonnegative.
					\end{Example}

\begin{Example}\label{exa:trid}
		Consider the    tridiagonal   matrix
		\be\label{eq:trida}
							A=\begin{bmatrix} 
																	a_1 & b_1 & 0 & \dots & 0											 \\
																	c_1 & a_2 & \ddots & \dots & \vdots											 \\
																	0 & \ddots & \ddots & \dots & \vdots											 \\
																	\vdots & \ddots & \ddots & \dots & b_{n-1}											 \\
																	0 & \dots & \dots & c_{n-1} & a_{n}											 
									\end{bmatrix} 
		\ee
		with~$ b_i,c_i\geq0$  for all~$i$.
		In this case, the dominance condition
		\be\label{eq:domi}
		a_i \geq b_i+c_{i-1} \quad \text{for all } i\in\{1,\dots,n\},
		\ee
		with~$c_0:=0$ and~$b_n:=0$, guarantees that~$A$ is TN~\cite[Ch.~0]{total_book}.
		\end{Example}

		An important subclass of TN matrices are the \emph{oscillatory matrices}  studied in the pioneering work
		of  \cite{gk_book}.
		A matrix~$A\in\R^{n\times n}$ is called \emph{oscillatory} if   $A$ is~TN and there exists an integer~$k>0$ such that~$A^k$ is~TP. It is well-known that a~TN matrix~$A$ is oscillatory if and only if it 
		is non-singular and irreducible~\cite[Ch.~2]{total_book}, and that in this case~$A^{n-1}$ is~TP. 
		\begin{Example}
		Consider  the matrix
		$
		A=\begin{bmatrix}   2&1&0 \\ 1&2&1\\0 &1 &2   \end{bmatrix}.
		$
		This matrix is~TN, non-singular and irreducible,
		so it is an oscillatory matrix. Here
	$
		A^{n-1}=A^2=\begin{bmatrix}   5&4&1 \\ 4&6&4 \\1 &4 &5   \end{bmatrix},
		$
		and it is straightforward to verify that this matrix is indeed~TP. 
		\end{Example}

 More generally,   the matrix~$A$ in~\eqref{eq:trida} satisfying~$b_i,c_i>0$, $i=1,\dots,n-1$,
and the dominance condition~\eqref{eq:domi}
	is~TN and irreducible. If it is also non-singular then it  is oscillatory.	

		An important property, that will be used throughout, is that
			the  product of two~TN [TP] matrices is a~TN [TP]  matrix. This   follows immediately from Def.~\ref{def:tntp}
			and the
		Cauchy-Binet formula for the minors of the product of two matrices~\cite[Ch.~0]{matrx_ana}. 
		(A straightforward proof of several important determinantal identities   
		can be found in~\cite{BRUALDI1983769}.)
		
		When using TN matrices to study dynamical systems, it is important to bear in mind that  in general 
		coordinate transformations do not preserve~TN (see Example~\ref{exa:scale} above). 
		An important exception, however, is positive diagonal scaling. Indeed, if~$D $ is a diagonal matrix with positive
		entries on the diagonal then multiplying a matrix~$A$ by~$D$ either on the left or
right   changes the sign of no minor, so in particular~$DAD^{-1}$ is TN [TP] if and only if~$A$ is~TN [TP].

		At this point we can already provide an intuitive explanation revealing  the so far unknown connection  
  between Smillie's results and~TN matrices. To do this, consider for simplicity the system~$\dot z=Jz$, with~$J$ a \emph{constant}
	    tridiagonal matrix
		  with positive entries on the super- and sub-diagonals. Then to first order in~$t$ we have~$z(t)=(I+tJ)z(0)$,
			and 
		  for any sufficiently
		small~$t>0$ 
		the matrix~$I+t J$ is~TN (see Example~\ref{exa:trid}), irreducible, and non-singular, so it is an oscillatory matrix.

		\subsection{Spectral properties of TN matrices}
		TN   matrices have a strong  spectral
structure: 
all their eigenvalues are real and nonnegative, and the corresponding eigenvectors
have special sign patterns. 
 This spectral structure is particularly evident in the case of oscillatory matrices.
\begin{Theorem}\label{thm:spec}\citep{Pinkus1996}
If~$A\in\R^{n\times n}$ is an oscillatory matrix then its eigenvalues are all real, positive, and distinct.
If we order the eigenvalues as~$\lambda_1>\lambda_2>\dots >\lambda_n>0$ 
and let~$u^k\in\R^n $ denote the eigenvector corresponding to~$\lambda_k$ then
for any~$1\leq i \leq j \leq n$ and any scalars~$c_i,\dots,c_j  $ that are not all zero
\[
			i-1 \leq s^-( \sum_{k=i}^j c_k u^k ) \leq s^+ (\sum_{k=i}^j c_k u^k )  \leq   j-1.
\] 
In particular,
\be\label{eq:suspu}
			s^-(u^i)=s^+(u^i)= i-1, \quad i=1,\dots,n.
\ee 
\end{Theorem} 
		
		Note that   eigenvectors may have some zero coordinates, but~\eqref{eq:suspu} implies that~$u^i \in \V$ for all~$i$. 
		\begin{Example}\label{exa:tpandspec}
		Consider the TP matrix
		$A =\begin{bmatrix}   5&4&1 \\ 4&6&4 \\1 &4 &5   \end{bmatrix} 
		$. Its eigenvalues are~$\lambda_1=2(3+2\sqrt{2})$, $\lambda_2=4$, $\lambda_3=2(3-2\sqrt{2})$ 		
		with corresponding eigenvectors 
	$u^1=\begin{bmatrix} 1 & \sqrt{2} & 1\end{bmatrix}' $, $u^2=\begin{bmatrix} -1 & 0 & 1\end{bmatrix}' $,
	and $u^3=\begin{bmatrix} 1 & -\sqrt{2} & 1\end{bmatrix}' $. Note that~$s^-(u^k)=s^+(u^k)= k-1$ for~$k=1,2,3$. 
		\end{Example}

		Note that~$A$ oscillatory implies  that~$A$ is componentwise nonnegative
		and that~$A^\ell$ is TP for some integer~$\ell \geq 1$ and thus by the Perron-Frobenius Thm.~\citep{matrx_ana} $A^\ell$
		(and thus~$A$) has a real positive eigenvalue~$\lambda_1$ with maximum modulus  and the corresponding eigenvector~$u^1$ has positive 
		entries, i.e.~$s^-(u^1)=s^+(u^1)=0$. 
		Thm.~\ref{thm:spec} shows that for oscillatory matrices a much stronger spectral structure arises.

		\subsection{Sign variation diminishing property of TN matrices  }
				TN   matrices enjoy a remarkable variety of mathematical properties. 
For our purposes, the most relevant property is that multiplication by a TN matrix cannot increase the 
sign variation of a vector. This is  
    the \emph{sign variation diminishing property}~(SVDP) of linear TN   transformations.

		Let~$A$ be a TN  EB matrix. 
		Pick~$x\in \R^n$, and let~$y:=Ax$. Then there exists at most one index~$i$ such  
		  that~$\sgn(y_i) \not = \sgn(x_i)$, and
			either~$y_i=x_i+p x_{i-1}$ or~$y_i=x_i+p x_{i+1}$, and since~$p\geq 0$
		the sign can change only in the ``direction'' of~$x_{i-1}$ or~$x_{i+1}$. 
		In either case, neither~$s^-$ or~$s^+$ may increase. 
		We  conclude that if~$A$ is TN  EB  then
		\begin{align}\label{eq:sminchange}
		s^-(Ax)&\leq s^-(x) \text{ for all }x\in\R^n,
		\end{align}
		and
		\begin{align}\label{eq:splushcnage}
		s^+(Ax)&\leq s^+(x) \text{ for all }x\in\R^n.
		\end{align}
		A similar argument shows that if~$A$ is TN GEB then~\eqref{eq:sminchange} holds. 
		However,~\eqref{eq:splushcnage} does not hold in general for a TN GEB matrix~$A$. 
		For example,~$A=0$
		 is   TN GEB and clearly~$s^{+}(Ax)=s^{+}(0)=n-1$ may be larger than~$s^{+}(x)$.
 		
			This SVDP  can be extended to all~TN matrices using the 
			following fundamental decomposition result.
   \begin{Theorem}~\cite[Ch.~2]{total_book}\label{thm:bifact}
		Any~TN matrix  can be expressed as a product of TN GEB matrices.
		\end{Theorem}
		\begin{Example}
		For the simplest example of this bidiagonal factorization, consider the case~$n=2$. Suppose that 
		$
			A=\begin{bmatrix} a_{11}& a_{12} \\a_{21} &a_{22}\end{bmatrix}
		$
		is~TN, that is,~$a_{ij}\geq 0$, $i,j=1,2$, and~$\det(A)=   a_{11}a_{22}-a_{12}a_{21}  \geq 0$. We consider two cases.  
		If~$a_{11}>0$ then
		$
		A=\begin{bmatrix}1&0 \\ a_{21}/a_{11} & 1 \end{bmatrix} 
		\begin{bmatrix}  a_{11}& 0 \\0& \det(A)/a_{11} \end{bmatrix}
			\begin{bmatrix}1& a_{12}/a_{11}  \\0  &1 \end{bmatrix} ,
		$
		and this is a product of the required form. If~$a_{11}=0$ then the TN of~$A$ implies that~$a_{12}a_{21}=0$. 
		Assume without loss of generality that~$a_{12}=0$. Then
		$
A=\begin{bmatrix} 0 &0 \\ a_{21}  & 1 \end{bmatrix} 
		\begin{bmatrix}  1& 0 \\0& a_{22} \end{bmatrix},
			$
and this is again 	a product of the required form.
\end{Example}

		Combining Thm.~\ref{thm:bifact} with the fact that~\eqref{eq:sminchange} holds for all TN GEB matrices 
				implies that~\eqref{eq:sminchange} holds for \emph{any}~TN matrix~$A$. 
		\begin{Remark}\label{rem:cnbs} 		
		As noted above,~$A$    TN does not imply~\eqref{eq:splushcnage}.
					However,~\eqref{eq:splushcnage}
					does hold   if~$A$ is a  \emph{nonsingular}~TN matrix~\cite[Ch.~4]{total_book}.  
		\end{Remark}
				
 TP matrices satisfy a  stronger 
SVDP.
\begin{Theorem}~\cite[Ch.~4]{total_book}\label{thm:tpscdp}
					If~$A \in \R^{n\times n}$ is TP then 
					\be\label{eq:strong}
									s^+(Ax)\leq s^-(x) \text{ for  all }  x\in \R^n \setminus\{0\}.
					\ee
					If~$A$ is TN and nonsingular then~$s^+(Ax)\leq s^-(x)$      
					for all~$x\in \R^n$ such that either~$x$ has no zero entries or~$Ax$ has no zero entries. 
\end{Theorem}


A natural question is whether   SVDP     characterizes~TN or~TP
matrices. 
Recall that a   matrix is called \emph{strictly  sign-regular}~(SSR) if for every~$k$
all
  minors of order~$k$ 
   are non-zero and share a common sign (that may vary from size to size). 
For example,
$
A=\begin{bmatrix}  1 & 2 \\ 3& 1   \end{bmatrix}
$ is~SSR because all 
  minors of order one are positive, and the single  minor of order two 
	is non-zero. 
Obviously,~TP   matrices are~SSR.
\begin{Theorem}\label{thm:regu}
Let~$A\in\R^{n \times n}$ be a nonsingular matrix. 
Then the following properties are equivalent. 
\begin{enumerate}[(a)]
							\item $A$ is SSR. \label{it:ssr}
							\item $s^+(Ax)\leq s^-(x)$ for all~$x\in \R^n\setminus\{0\}$. \label{it:svdpstrong}
\end{enumerate}
\end{Theorem}
Note that  the first assertion  in Thm.~\ref{thm:tpscdp} follows from this result.
For the sake of completeness, we detail the proof of Thm.~\ref{thm:regu}
   in the Appendix. 
The assumption that~$A$ is nonsingular cannot be dropped. 
The  next example demonstrates this.
\begin{Example}
Consider  the matrix
$
A=\begin{bmatrix}  2 & 2 \\ 1& 1   \end{bmatrix}
$.  
 Pick~$x\in\R^2\setminus \{0\}$.
If~$s^+(Ax)=0$   then clearly~$s^+(Ax)\leq s^-(x)$.
Otherwise,~$s^+(Ax)=1$, that is, either~$2x_1+2x_2 \leq 0$ and~$x_1+x_2 \geq 0$ 
or~$2x_1+2x_2 \geq 0$ and~$x_1+x_2 \leq 0$.
Both these  cases imply that~$x_1=-x_2$, and since~$x\not =0$, we conclude that~$s^-(x)=1$.
Thus, condition~\eqref{it:svdpstrong}  holds, but
condition~\eqref{it:ssr} does not hold,   as~$\det(A)=0$.
\end{Example}

\subsection{Dynamics of compound matrices} 
  Consider the matrix differential  equation~$\dot Y(t)=A(t)Y(t)$. What is the   dynamics of some  minor  of~$Y$? It turns out that we can express the dynamics of every~$p\times p$
minor in terms of all the~$p\times p$ minors of~$Y$ and the~$n^2$ entries of~$A$. 
To explain this we  
 review multiplicative and additive compound matrices and their role in certain differential equations~\citep{muldo1990}.

Given~$A\in\R^{n\times n}$ and~$p\in\{1,\dots,n\}$, consider the
$\binom{n}{p}^2$
 minors of~$A$ of size~$p\times p$. 
Each minor is defined by a set of row indexes~$1\leq i_1<i_2<\dots<i_p\leq n$ and column indexes~$1\leq j_1<j_2<\dots<j_p\leq n$. This minor 
is denoted by~$A(\alpha|\beta)$ where~$\alpha:=\{i_1,\dots,i_p\}$ and~$\beta:=\{j_1,\dots,j_p\}$.

The~$p$th \emph{multiplicative  compound matrix}~$A^{(p)}$ is the~$\binom{n}{p}\times  \binom{n}{p}$ matrix
that 
includes all these minors ordered lexicographically. For example, for~$n=3$ and~$p=2$ there are nine minors. 
The~$(1,1)$ entry of~$A^{(2)}$ is~$A(\{1,2\}|\{1,2\})$,
the~$(1,2)$ entry of~$A^{(2)}$ is~$A(\{1,2\}|\{1,3\})$,
and entry~$(3,3)$ of~$A^{(2)}$
  is~$A(\{2,3\}|\{2,3\})$.

An important property that follows from the Cauchy-Binet formula is 
\be\label{eq:pouy}
(AB)^{(p)}=A^{(p)} B^{(p)}.
\ee
This justifies the term multiplicative compound. 

The $p$th \emph{additive compound matrix} of~$A$
is   defined  by
$
				A^{[p]}:= \frac{d}{dh}  (I+hA)^{(p)} |_{h=0}.
$
Note that this implies that
\be\label{eq:poyrt}
(I+hA)^{(p)}= I+h  A^{[p]} +o(h)  ,
\ee
 and that~$A^{[p]}= \frac{d}{dh} \left(  \exp(h A)  \right)^{(p)} |_{h=0}$.

\begin{Example}
Consider the case~$n=3$ and~$p=2$. Then~$(I+hA)^{(2)}$ is the matrix depicted on top of p.~\pageref{eq:long1}
\begin{figure*}\label{eq:long1}
\[
			\kbordermatrix{ & 12 & 13 & 23  \\
												12&			(1+ha_{11})(1+h a_{22})-h^2 a_{21}a_{12} &  (1+ha_{11})ha_{23}-h^2 a_{13}a_{21}
												& h^2 a_{12}a_{23}-ha_{13}(1+ha_{22})\\
	13&	(1+ha_{11}) h a_{32}-h^2 a_{12}a_{31} &  (1+ha_{11}) (1+ha_{33})-h^2 a_{13}a_{31} 
	     & h a_{12} (1+ha_{33}) -h^2 a_{13}a_{32}\\
	23&		h^2  a_{21}a_{32}-h a_{31} (1+ha_{22})  &  h a_{21} (1+ha_{33}) -h^2 a_{23}a_{31} & (1+ha_{22})(1+ha_{33})-h^2a_{23}a_{32}
		} 
\]
\hrule
\end{figure*}
where the row [column] marks  are the indexes in~$\alpha$ [$\beta$].
Thus,
$A^{[2]}=\begin{bmatrix}
														 a_{11} + a_{22}  &    a_{23} &  -a_{13}  \\
		  a_{32} &  a_{11} + a_{33}  & a_{12}   \\
					  -a_{31}  &   a_{21}    &  a_{22}+a_{33} 
		\end{bmatrix}.
$
\end{Example}

Applying Cauchy-Binet  again gives 
$
							(I+hA)^{(p)} (I+h B)^{(p)}	=     ( I + hA+hB  + o(h) ) ^{(p)}  ,		
$
and this yields
\[
(A+B)^{[p]}= A^{[p]}+B^{[p]},
\]
justifying the term additive compound.

The additive compound is important when studying the dynamics of the multiplicative compound. 
For a time-varying matrix~$Y(t)$ we use the notation~$Y^{(p)} (t)$ for~$(Y(t))^{(p)} $. 
Suppose that \updt{$\dot Y(t)=A(t)Y(t)$}. Then 
\begin{align*}
      Y^{(p)} (t+h)&=  (Y(t)+h A(t)Y(t)   ) ^{(p)}+\updt{o(h) }\\
			&=(I+h A(t)   ) ^{(p)} Y  ^{(p)}(t)+\updt{ o(h) },
\end{align*}
and combining this with~\eqref{eq:poyrt} gives
\be\label{eq:povt}
\frac{d}{dt} Y^{(p)}(t)=A^{[p]}(t) Y^{(p)}(t).
\ee
Thus, the dynamics of~$Y^{(p)}(t)$ is also linear,  
  with the dynamical  matrix~$A^{[p]}(t) $.
	Note that combining~\eqref{eq:povt} and~\eqref{eq:pouy}
	implies that for a constant matrix~$A$,
	\[
	\exp(A^{[p]})=(\exp(A))^{(p)}.
	\]

\updt{
The matrix~$A^{[p]}$ can be determined explicitly. 
\begin{Lemma}\label{lem:poltr}
The entry of~$A^{[p]}$ corresponding to~$(\alpha|\beta)=(i_1,\dots,i_p|j_1,\dots,j_p) $  is:
\begin{itemize}
\item $\sum_{k=1}^p a_{i_k i_k}$ if ~$i_\ell=j_\ell$ for all~$\ell=1,\dots,p$; 
\item $(-1)^{\ell+m} a_{i_\ell j_m}$ 	 if all the indexes in~$\alpha$ and~$\beta$ coincide except for a single index~$i_\ell\not = j_m $; and
 \item $0$ otherwise. 
\end{itemize}
\end{Lemma} 
}
Note that the first case here  corresponds to diagonal entries of~$A^{[p]}$.  
Lemma~\ref{lem:poltr} is  usually proven by manipulating determinants~\citep{schwarz1970}
or using exterior powers~\citep{fiedler_book}.

\begin{Example}
For~$p=1$  we have~$Y^{(1)}=Y$ and Lemma~\ref{lem:poltr}
yields~$A^{[1]} = A$, so we
obtain~$\dot Y=AY$.  
\end{Example} 
\begin{Example} 
For~$p=n$ the matrix~$Y^{(n)}$ includes a single entry which is~$\det(Y)$, whereas
Lemma~\ref{lem:poltr} yields~$A^{[n]}=\trace(A)$,
so~\eqref{eq:povt}  yields  
\be\label{eq:ajl}
\frac{d}{dt} \det(Y(t))= \trace(A(t)) \det( Y(t) ),
\ee	 
which is the  Abel-Jacobi-Liouville identity.
\end{Example}

   \cite{schwarz1970} considered the following problem. Suppose that~$\dot Y=AY$. 
Then the dynamics of the~$p$th multiplicative compound~$Y^{(p)}$ is given by the linear system~\eqref{eq:povt}. 
When will \emph{every}~$A^{[p]}$ be a Metlzer matrix?  An interesting property that will be proven below is 
that  if~$A^{[1]}$ and~$A^{[2]}$ are Metzler then every~$A^{[p]}$ is Metzler.  
The next example demonstrates this for the case~$n=4$. 

\begin{Example}
 Consider the case~$n=4$ i.e~$A=\{a_{ij}\}_{i,j=1}^4$. 
Then Lemma~\ref{lem:poltr} yields
\begin{align*}
A^{[2]}&=\left [ \begin{smallmatrix} 
																 a_{11}+a_{22}   & a_{23} & a_{24} & -a_{13} & -a_{14} & 0 \\
																a_{32}      &a_{11}+a_{33} & a_{34} & a_{12} & 0 & -a_{14} \\   
																a_{42}     &a_{43}   &a_{11}+a_{44} & 0 & a_{12} & a_{13} \\   
																-a_{31}      &a_{21} & 0 & a_{22}+a_{33} & a_{34} &  -a_{24} \\
																-a_{41}      &0 & a_{21} & a_{43} & a_{22}+a_{44} &  a_{23} \\ 
																0      &-a_{41} & a_{31} & -a_{42} & a_{32} &  a_{33}+a_{44}  
  \end{smallmatrix}\right] , \\
A^{[3]}&=\left[ \begin{smallmatrix}     
																a_{11}+a_{22}+a_{33} & a_{34} & -a_{24} & a_{14} \\
																a_{43}& a_{11}+a_{22}+a_{44} & a_{23} & -a_{13}  \\
																-a_{42}& a_{32} & a_{11}+a_{33}+a_{44} & a_{12}  \\
                                a_{41}& -a_{31}& a_{21}& a_{22}+a_{33}+a_{44} 
 \end{smallmatrix} \right].
\end{align*}
Suppose that~$A^{[1]}=A$ is Metzler, that is,~$a_{ij}\geq 0$ for all~$i\not =  j$. 
Then~$A^{[2]}$ is Metzler iff~$a_{13}=a_{14}=a_{24}=a_{31}=a_{41}=a_{42}=0$. Under these conditions, we see that~$A^{[3]}$ 
is also Metzler. 
The matrix~$A^{[4]}$ is also  Metzler, as it is a scalar. 
\end{Example}

We note in passing that contractivity  of  
 the second additive compound~$J^{[2]}$, where~$J$ is the Jacobian  in the variational equation~\eqref{eq:zdoteqjj},
can be used to analyze the existence and stability of  (non-constant) periodic solutions of~\eqref{eq:fsys}~\citep{muldo1990}. This
proved useful in  many models     from systems biology, see e.g.~\citep{tri_feed,SEIR_LI_MULD1995,mosquito_2017}.
		
	\section{Totally positive differential systems}\label{sec:main}

Consider  the matrix system
\be\label{eq:yata}
\dot Y(t)=A(t)Y(t),\quad Y(t_0)=I, 
\ee
with~$A(t)$ a continuous function of~$t$.
Let~$Y(t,t_0)$ denote 
 the solution of~\eqref{eq:yata} at time~$t$.
 \cite{schwarz1970}
called~\eqref{eq:yata}  a  totally  nonnegative  differential system~(TNDS)
if for every~$t_0$ the solution~$Y(t,t_0)$ is~TN   for  all~$t \geq t_0$, and a
   totally  positive   differential system~(TPDS)
if for every~$t_0$ the solution~$Y(t,t_0)$ is~TP   for  all~$t >  t_0$.

Schwarz  combined the   Peano-Baker representation
for the solution of~\eqref{eq:povt} (see, e.g.~\cite{linear_st}) 
 and Lemma~\ref{lem:poltr}  
to derive necessary and sufficient conditions for a system to be~TNDS [TPDS].
Stated in modern terms, his analysis  is based on the fact that~\eqref{eq:yata}
is TNDS [TPDS] iff~\eqref{eq:povt} is a 
 cooperative [strongly cooperative] dynamical 
 system for all~$p$. 
Indeed, note that~$Y(t_0)=I$ implies that~$Y^{(p)}(t_0)=I$ for all~$p$, so in particular all the minors at time~$t_0$ are nonnegative.
The cooperative [strongly coopertive] dynamics maps these nonnegative initial conditions at time~$t_0$ 
 to a nonnegative [positive] solution for all~$t>t_0$. 
Schwarz also studied the implications of TNDS/TPDS of~\eqref{eq:yata} on the number of sign variations
in a   solution of the associated vector differential equation.

Fix  a time interval~$(a,b)$  with~$-\infty\leq a<b \leq \infty$. 
For any pair~$t_0,t$, with~$a< t_0 \leq t  <b$, consider the vector differential equation
\be\label{eq:azz}
\dot z(t)=A(t)z(t), \quad z(t_0)=z_0. 
\ee
We  assume throughout a more general case than in~\cite{schwarz1970}, namely, that  
\begin{align}\label{eq:sas}
\text{ }  A:(a,b)& \to \R^{n\times n} 
\text{ is a  matrix of locally (essentially) } \nonumber \\&\text{bounded  measurable 
functions.}
\end{align}
This more general formulation is important in the context of control systems. Indeed,   consider~$\dot x=f(x,u)$, with~$u$ a control input. Then \updt{$z:=\dot x$ satisfies the} variational equation~$\dot z=J(x,u)z$, where~$J$ is the Jacobian of~$f$
with respect to~$x$. In many cases, for example when considering optimal controls, one must allow measurable controls (see, e.g.~\cite{liberzon_opt_cont})  
and thus~$t\to J(t)$ is typically a measurable, but not continuous, matrix function.

\updt{The transition matrix associated with~\eqref{eq:azz} is defined by}   
\be\label{eq:defuder}
					\frac{d}{dt} \Phi(t,t_0)=A(t) \Phi(t,t_0), \quad \Phi(t_0,t_0)=I. 
\ee
 Recall that~\eqref{eq:sas} implies that~\eqref{eq:defuder} admits a unique, locally absolutely continuous,    nonsingular 
solution for   all pairs~$ (t_0,t) \in(a,b)\times(a,b)$ 
(see, e.g.,~\cite[Appendix~C]{sontag_book}). 

The formula~$z(t)= \Phi(t,t_0)  z(t_0)$ suggests
that if~$\Phi(t,t_0)$ is~TP then~$\sigma(z(t ))$ will be no larger than~$\sigma(z(t_0 ))$.
  The next result formulates this idea.
	\begin{Theorem}\label{thm:exp_is_tp}
Consider the time-varying linear system:
\be\label{eq:linti}
\dot z(t)=A(t)z(t), 
\ee
with~$A(t)$ satisfying~\eqref{eq:sas} 
and suppose that
\be \label{eq:phist}
\Phi(t,t_0) \text{ is TP for all } a<  t_0 <t<b   . 
\ee
If~$z(t)$ is not the trivial solution~$z(t)\equiv  0$  then:
\begin{enumerate}[(1)]
\item the functions~$s^-(z(t)) , s^+(z(t))$ are   non-increasing functions of
 time on~$(a,b)$;
\item $z(t) \in \V$  for all~$t\in (a,b)$, except perhaps for up  to~$n-1$ discrete  values of~$t$. 
\end{enumerate}
\end{Theorem}

As we will see in Section~\ref{sec:app} below, both  these    properties  are useful in the analysis       
of nonlinear~ODEs. 

{\sl Proof of Thm.~\ref{thm:exp_is_tp}.}
  For any~$a<t_0<t<b$ we have
  $z(t) =\Phi(t,t_0)z(t_0 )$ and since 
	the matrix here is~TP,~\eqref{eq:strong} yields
	\be\label{eq:slop}
	s^+(z(t))\leq s^- (z(t_0) ) ,
	\ee
	and thus~$s^-(z(t))\leq s^+(z(t))\leq s^- (z(t_0) ) $.
	If~$z( t_0) \in \V$ then~$s^- (z(t_0) )  = s^+ (z(t_0) )$, so~\eqref{eq:slop} yields  
	\[
					s^+(z(t))  \leq  s^+ (z(t_0) ) .
	\]
	If~$z(t_0)\not \in \V$ then~$s^- (z(t_0) )  < s^+ (z(t_0 ) )$, so 
	\[
					s^+(z(t))  < s^+ (z(t_0) ). 
	\]
	Thus, $s^+(z(t))$ never increases, 
	and it strictly decreases as~$z(t)$ goes through a point that is not in~$\V$. 
	Since~$s^+$ takes  values  in~$\{0,1,\dots,n-1\}$, this implies
	that~$z(t)\in\V$ for all~$t $, except perhaps for up to~$n-1$ discrete  points. \hfill{$\square$}

The  proof  of Thm.~\ref{thm:exp_is_tp}
 shows that we may view~$s^- (z(t))$ and~$s^+(z(t))$ as   integer-valued Lyapunov functions of
 the time-varying linear system~\eqref{eq:linti}. The same is true for~$\sigma(z(t))$.
Indeed, if~$z(\tau)\not \in\V$ for some~$\tau$ (and recall that this can only hold for up to~$n-1$ discrete  points)
then~$z(\tau^-),z(\tau^+) \in \V$, so
\[
				\sigma(z(\tau^+))=s^+(z(\tau^+)) < s^+(z(\tau^-))= \sigma(z(\tau^-)).
\]
Thus,~$\sigma(z(t))$ is piecewise constant with no more than~$n-1$ points of discontinuity and at these points it strictly decreases.

\begin{Remark}\label{rem:conv}
Using Thm.~\ref{thm:regu} yields  a   converse for
Thm.~\ref{thm:exp_is_tp}.
 Indeed, suppose that  the solution of~$\dot z =A z$   satisfies~\eqref{eq:slop} for 
 all~$a<t_0<t<b$ and all~$z(t_0)\in\R^n \setminus\{0\}$. Then  
	using the fact that~$z(t)=\Phi(t,t_0)z(t_0)$ and Thm.~\ref{thm:regu} 
imply that~$\Phi(t,t_0)$ is~SSR for all~$a<  t_0 < t<b$. 
Pick~$1\leq p\leq n$. 
Then all the~$\binom{n}{p}$ minors of order~$p$ of~$\Phi(t,t_0)$ are either all positive or all negative. 
Since the matrix~$\Phi(t_0,t_0)=I$ has a minor of order~$p$ that is one,
and the~$SSR$ property means that for any~$t>t_0$ this  minor   is not  zero, we conclude by continuity 
that this  minor  of~$\Phi(t,t_0)$   is     positive for all~$t>t_0$, and thus all minors of order~$p$ are positive
for all~$t>t_0$.
  Since this holds for all~$p$,~$\Phi(t,t_0)$ is~TP  for  all~$a<t_0<t<b$. 
\end{Remark}
In particular, 
 the next example shows that if we change the word~``TP'' in condition~\eqref{eq:phist}
 to ``TN'' then   Thm.~\ref{thm:exp_is_tp} no longer 
holds.

\begin{Example}
Consider the constant matrix
$
			A = \begin{bmatrix} a_{11} & a_{12} \\0 &a_{22} \end{bmatrix},
$
with~$a_{12} >  0$.   For $t_0=0$,
\[
	\Phi(t,t_0)=		\exp(At)=\begin{bmatrix}1+ a_{11}t+o(t) & a_{12}t+o(t)  \\0 &1+a_{22}t+o(t) \end{bmatrix} ,
\]
and thus there exists~$T>0$ such  
  that~$\exp(At)$ is~TN for all~$t\in[0,T]$.
	However,~$\exp(At)$ is not~TP for any~$t$. 
For~$z(0)=\begin{bmatrix} 1&0 \end{bmatrix}'$,
the solution of~$\dot z= A z$ is~$z(t)= \exp(a_{11}t) \begin{bmatrix} 1&0 \end{bmatrix}'$.
Since~$z_2(t)\equiv 0$,~$z(t) \not \in \V$ for all~$t\geq 0$.
\end{Example}

We now formally state the definitions of a~TNDS and a TPDS. 
	\begin{Definition}
	We say   that~\eqref{eq:defuder} is a~TNDS if for all~$ a  <t_0\leq t <b$ 
	the matrix~$\Phi(t,t_0)$ is~TN. 
	We say   that~\eqref{eq:defuder} is a~TPDS if for all~$ a<  t_0   < t<b $ 
	the matrix~$\Phi(t,t_0)$ is~TP. 
\end{Definition}

\begin{Example}\label{exa:cosh}
							Consider  the matrix~$A(t)=\begin{bmatrix}0&t\\t&0\end{bmatrix}$.
								Note that this is tridiagonal and with positive entries on the sub- and super-diagonals 
								for all~$t>0$.
							  In this case, the solution of~\eqref{eq:defuder}  is
							\[
							\Phi(t,t_0)=\begin{bmatrix}
							\cosh( (t^2-t_0^2)/2 ) & \sinh( (t^2-t_0^2)/2 ) \\
							\sinh( (t^2-t_0^2)/2 ) & \cosh( (t^2-t_0^2)/2 ) 
							\end{bmatrix}.
							\]
Note that every entry here is positive for all~$t>t_0\geq 0 $ and that~$\det(\Phi(t,t_0))\equiv 1$, so~$\Phi(t,t_0)$ is~TP on any interval~$(a,b)$, with~$a\geq 0$.
Thus the system is a~TPDS on such an interval. 

For~$A(t)=\begin{bmatrix}0&t\\0&0\end{bmatrix}$  
							  the solution of~\eqref{eq:defuder}  is
							$
							\Phi(t,t_0)=\begin{bmatrix}
							 1 & (t^2-t_0^2)/2  \\
							  0& 1 
							\end{bmatrix}
							$. This matrix is~$TN$ (but not~TP) for all~$t \geq t_0\geq 0 $,
							so the system is a~TNDS on any interval~$(a,b)$, with~$a\geq 0$. 
\end{Example}

 Since the
 product of~TN [TP] matrices is a~TN [TP] matrix,
a sufficient condition for TNDS  [TPDS] is that   there exists~$\delta>0$, that does not depend on~$t_0$,
 such that for any~$\epsilon \in(0,\delta)$
\[
				\Phi(t_0+\varepsilon,t_0) \text{ is TN [TP]  for all }
				a< t_0 <   b -\varepsilon .
\]
 
Our next goal is to describe    conditions on~$A(t)$
guaranteeing that~\eqref{eq:defuder} is a~TNDS or a~TPDS.
 A related question has already been addressed  by
\cite{Loewner1955} who studied the infinitesimal generators of the group of~TN matrices.


It is useful to first  consider      the case of a constant matrix.  We require the following notation.
\begin{Definition}
Let~$\M \subset \R^{n\times n}$ [$\M^+ \subset \R^{n\times n}$] denote the set of tridiagonal matrices
with nonnegative [positive] entries on the sub- and super-diagonal. 
\end{Definition}

The next result provides a simple necessary and sufficient condition
 for~\eqref{eq:defuder}, with~$A$ a \emph{constant} matrix,
to be~TNDS or~TPDS.
 \begin{Theorem}\label{thm:tic}\citep{schwarz1970}
Fix an interval~$(a,b)$.
The system~$\dot U(t)=AU(t)$ is TNDS [TPDS] on~$(a,b)$
if and only if~$A\in \M$
[$A\in  \M^+$]. 
\end{Theorem} 

Due to the importance of this result, we provide two different     proofs.  
The second proof follows~\citep{schwarz1970}
and is   useful when we consider below  the case where~$A$ is time-varying. 
The first proof is  new  and possibly     easier to follow.

{\sl First Proof of Thm.~\ref{thm:tic}.}
The solution of the   matrix differential equation:
\be\label{eq:vw}
			\dot U(t)=A U(t), \quad U(0)=I, 
\ee
is
\be\label{eq:uuexp}
U(t) =I+\frac{A t}{1!}+\frac{A^2 t ^2}{2!}+\dots.
\ee
If~$a_{ij}<  0$ for some~$i\not =j$ then~$u_{ij}(t)<0$ for 
all~$t>0$ sufficiently  small. Thus, a necessary condition for~$U(t)$ to be~TN
for all~$t>0$ sufficiently small is that~$A$
is a Metzler matrix. 
  
We now show that another  necessary condition for~$U(t)$ to be~TN
for all~$t>0$ sufficiently small is that~$A$
is tridiagonal.  If~$n=2$ then~$A$ is always tridiagonal, so   assume that~$n\geq 3$. 
Pick~$1\leq j<k<i\leq n$. Then~\eqref{eq:uuexp} yields 
\begin{align*}
								  \det   \begin{bmatrix} u_{kj}(t) &u_{kk}(t)\\u_{ij}(t)&u_{ik}(t)\end{bmatrix}   
								&=t^2 a_{kj}  a_{ik}-(1+t a_{kk})t a_{ij}+o(t)\\
								&=-t a_{ij}+o(t).
\end{align*}
This implies that if~$a_{ij}>  0$ for some~$i>j+1$ then~$U(t)$ has 
a negative $2\times 2 $  minor    for all~$t >0$ sufficiently  small.
A similar argument shows that 
if~$a_{ij}>  0$ for some~$j>i+1$ 
then~$U(t)$ has 
a negative~$2\times 2 $ minor    for all~$t >0$ sufficiently  small.

Summarizing, a necessary condition for~$U(t)$ to be~TN for all~$t >0$ sufficiently  small 
is that~$A\in\M$. 
Suppose that this indeed holds. \updt{If~$a_{ij}=0$ for some~$i,j$ with~$|i-j|=1$ 
 then the tridiagonal structure implies that for any~$k$, entry~$i,j$ of~$A^k$ is also zero, 
so~$u_{ij}(t)=0$ for all~$t\geq 0$, and thus~$U(t)$ is not~TP.}
We conclude that a necessary condition for~$U(t) $ to be~TN  [TP] for all~$t>0$ 
is that~$A\in \M$ [$A\in  \M^+$]. 

To prove the converse implication, assume that~$A\in \M^+ $.
Then~\eqref{eq:uuexp} implies that~$U(t)$ is irreducible for all~$t>0$. 
 It is
enough to show that there exists~$\varepsilon_0>0$ such that~$U(t)$ is~TP 
for all~$t\in (0,\varepsilon_0)$. 
By the result stated in Example~\ref{exa:trid}, the matrix~$I+\frac{A t }{k}$ is~TN for any  sufficiently
large~$k$. 
Using the formula
$
				U(t)=\lim_{k\to \infty} \left(I+\frac{At}{k} \right )^k, 
$
we conclude that  there exists~$\varepsilon>0$ such that 
$U(t)$ is~TN for all~$t \in  \I:=(0,\varepsilon )$.
 Summarizing,~$U(t)$ is irreducible, non-singular and~TN
for all~$t \in \I  $ and is thus an oscillatory matrix on this time interval. 
Hence,~$(U(t))^{n-1}$ is~TP for all~$t\in \I $,
 so~$U(t) $ is~TP for all~$ t\in(0, (n-1)\varepsilon) $.
This completes the     proof  for the~TPDS case. The~TNDS case   follows similarly.~\hfill{$\square$}

\begin{Remark}
The proof above has an important and non-trivial implication.
				It shows that a necessary condition for~$U^{(1)}(t)$ and~$U^{(2)}(t)$ 
				to be componentwise nonnegative [positive] for all~$t>0$ sufficiently small 
				is that~$A\in \M$ [$A\in\M^+$].
				On the other hand, if~$A\in M$ [$A\in\M^+$] then \emph{every} minor of~$U(t)$ is
				nonnegative [positive] for all~$t>0$. Thus, checking the~$1\times 1$ and~$2\times 2$
				minors of~$U(t)$ is enough to establish~TNDS or TPDS. 
\end{Remark}

We now provide another  proof of Thm.~\ref{thm:tic} that
 follows the ideas in~\cite{schwarz1970}.

{\sl Second Proof of Thm.~\ref{thm:tic}.}
Given~\eqref{eq:vw},  recall that for any~$1\leq p\leq n$ 
  the induced dynamics for 
the~$p\times p$ minors of~$U(t)$ is
\be\label{eq:poxrr}
			\dot U^{(p)} =A^{[p]} U^{(p)}, \quad U^{(p)}(0)=I, 
\ee
where~$\dot U^{(p)}:=\frac{d}{dt} (U^{(p)}(t))$, and~$A^{[p]}$
is given in Lemma~\ref{lem:poltr}.

Our goal is to find conditions guaranteeing that~$A^{[p]}$, $p=1,\dots,n$,
 is Metzler. Indeed,
it is easy to see that if an   off-diagonal entry of~$A^{[p]}$ is negative
then there will be an entry in~$U^{(p)}$ that is negative for all~$t>0$ sufficiently small. If~$A$ is not Metzler then~$A^{[1]}=A$ has a negative off-diagonal entry, so we conclude that a necessary condition
 for~TNDS is that~$A$
is Metzler. 

Assume that~$a_{ij}>0$ for some~$i,j$
with~$|i-j|>1$. Consider~\eqref{eq:poxrr} with~$p=2$.
If~$j\geq i+2$ [$i\geq j+2$]  then Lemma~\ref{lem:poltr} yields that the off-diagonal entry of~$A^{[2]}$ corresponding to~$(i,i+1|i+1,j)$ [$(j+1,i|j,j+1)$]   is~$(-1)^{1+2} a_{ij}<0$ [$(-1)^{2+1} a_{ij}<0$]. We conclude that a necessary condition for~TNDS is that~$A \in \M$. 
If~$A \in \M$ and~$a_{ij}=0$
for some~$i,j$ with~$|i-j|=1$ then~$U_{ij}(t)\equiv 0$, so we conclude that
a  necessary condition for~TPDS is that~$A \in \M^+$.

To prove the converse implication, assume that~$A\in \M$. 
Pick~$1\leq p\leq n$. We now show that~$A^{[p]}$ 
 is  Metzler. By Lemma~\ref{lem:poltr}, an off-diagonal
entry of~$A^{[p]}$ corresponding to~$(\alpha|\beta)$
is either zero and then we are done, or
it is~$
					(-1)^{\ell+m} a_{i_\ell j_m}
$,
when~$\alpha$ and~$\beta$ have~$p-1$ identical entries 
 and~$i_\ell\not = j_m$.
Since~$A\in\M$, the term~$	(-1)^{\ell+m} a_{i_\ell j_m}$ can  be 
 negative   only if the set~$\{i_\ell,j_m\}$
is the set~$\{k,k+1\}$ for some~$k\in\{1,\dots,n-1\}$. 
Assume this is so. Then~$\alpha$ and~$\beta$ each include~$p$ increasing indexes, $p-1$ of these coincide, with one  of~$k,k+1$ appearing in~$\alpha$ but not in~$\beta$ and the second appearing in~$\beta$ but not in~$\alpha$. 
  This implies  that~$\ell=m$. But then~$(-1)^{\ell+m} 
a_{i_\ell j_m}=a_{i_\ell j_m}\geq 0$. 
We conclude that if~$A\in \M$ then~$A^{[p]}$ is Metzler and then it follows from   known results on
 cooperative dynamical systems (see, e.g.~\cite{hlsmith})
that every entry of~$Y^{(p)}(t)$ is nonnegative 
for all~$t\geq 0$. Since this holds for every~$p$, the system is~TNDS. 

Assume now that~$A\in \M^+$. 
Then~$A^{[1]}=A$ is an irreducible matrix. Pick~$2\leq p \leq n$. 
We will show that~$A^{[p]}$
is irreducible using the equivalence  between irreducibility
and strong connectivity 
 of the 
 adjacency graph associated with~$A^{[p]}$
(see e.g.~\cite[Ch.~6]{matrx_ana}). As we only verify strong connectivity of this  graph, it 
    is enough to consider  the case where~$A$ is tridiagonal, with the main diagonal all zeros, and the sub- and super-diagonal is all ones. In this case all the entries of~$A^{[p]}$ are zero, except for those that correspond to~$(\alpha,\beta)$ where
 exactly~$p-1$ entries of~$\alpha$ and~$\beta$   coincide, and the two remaining indexes are~$i_\ell$ and~$j_m=j_\ell$
 satisfy~$|i_\ell-j_\ell|=1$. 
Then~$A^{[p]} (\alpha|\beta)=1$. 

Consider the adjacency  graph associated with
the matrix~$A^{[p]}$. Every node in this graph corresponds to
 a set of~$p$  increasing indexes~$1\leq i_1< \dots<i_p\leq n$, and there are~$\binom{n}{p}$ nodes. There is a undirected edge 
between nodes~$\alpha$ and~$\beta$ if exactly~$p-1$ entries of~$\alpha$ and~$\beta$   coincide, and the two remaining indexes 
 satisfy~$|i_\ell-j_\ell|=1$. This means that there is a path in the graph from every node  to the node~$(1,2,\dots,p)$. Hence, the graph is strongly connected, so~$A^{[p]}$ is irreducible. 
We conclude that if~$A\in \M^+$ then~$A^{[p]}$ is Metzler and irreducible, and  it follows from   known results on
 cooperative dynamical systems (see, e.g.~\cite{hlsmith})
that every entry of~$U^{(p)}(t)$ is positive 
for all~$t> 0$. Since this holds for every~$p$, the system is~TPDS.~\hfill{$\square$}

\begin{Example}
Consider the case~$n=3$ and the matrix
$
A=\begin{bmatrix}     0& a_{12}& 0 \\ a_{21} &0& a_{23} \\ 0& a_{32} & 0   \end{bmatrix}  
$,
with~$a_{ij}>0$. Note that~$A\in \M^+$.
The  solution~$U(t)$ of~\eqref{eq:vw} with~$t_0=0$ is  
\begin{align*}
			 &I+At+ A^2 \frac{t^2 }{2}       +o(t^2) \\
			&=   \begin{bmatrix}    1+ a_{12} a_{21} \frac{t^2}{2}     & 
			                            a_{12}  t  &
																		a_{12} a_{23}\frac{t^2}{2}     \\
																		 a_{21}  t  &
																		  1+ (a_{12} a_{21} + a_{23}a_{32}  )\frac{t^2}{2}     &
																			a_{23} t \\
																			 a_{21} a_{32} \frac{t^2}{2}     &
																			a_{32} t &
																			1+ a_{23}  a_{32}\frac{t^2}{2}      
																		      \end{bmatrix}  \\&                                        +o(t^2). 
\end{align*}
It
 is straightforward to see that every~$1\times 1$ and~$2\times 2$
minor here is positive for all~$t>0$ sufficiently small. 
Also, the  Abel-Jacobi-Liouville identity~\eqref{eq:ajl}
yields~$\det(U(t)) \equiv \det(U(0)) = 1$,  
 so~$U(t)$ is~TP for all~$t>0$ sufficiently small. 
\end{Example}

We now turn   to
consider the time-varying case. Here we generalize the 
results in~\cite{schwarz1970} to our more general, measurable case. 

 \begin{Theorem}\label{thm:aprre}
 Fix an interval~$(a,b)$.
The system~\eqref{eq:defuder} with~$A(t)$ satisfying~\eqref{eq:sas} is TNDS
on~$(a,b)$ iff~$A(t) \in \M$ for  almost all~$t\in(a,b)$.
\end{Theorem}

	To prove this, we  require the following result. We use~$Q\geq 0$ [$Q\gg 0$] to denote that every 
	entry of the matrix~$Q$ is nonnegative [positive].
	\begin{Lemma}\label{lem:phiplis}
          For any $t_0$ and $t$ with $a\leq t_0 < t \leq b$, denote by
          $\Theta(t,t_0)$ the unique solution, at time $t$,
          of $\dot \Theta(s) =A(s)\Theta(s)$,~$\Theta(t_0)=I$. 
       Then the following two conditions are equivalent.
	\begin{enumerate}[(1)]
	 \item \label{enu:tnds} $\Theta(t,t_0)\geq 0$ for all~$a<t_0\leq t<b$; 
	 \item \label{enu:ametz} $A(t)$ is Metzler for almost all~$t \in(a,b)$. 
	\end{enumerate}
	\end{Lemma}

	{\sl Proof of Lemma~\ref{lem:phiplis}.}
Assume that~$A(t)$ is Metzler for almost all~$t \in(a,b)$.  Since $A(t)$ is a
matrix of locally (essentially) bounded measurable functions, we may pick an
$r>0$ such that~$rI+A(t)\geq0$ for almost all $t\in[a,b]$.
Pick~$t_0\in(a,b)$.
We begin by assuming that~$\Theta(t_0)\gg 0$,
and introduce the auxiliary matrix function
$\Psi(t):=e^{r(t-t_0)}\Theta(t,t_0)$. 
Suppose that there would exist some $t_1>t_0$ such that
$\Theta(t_1,t_0)\not\geq0$, or equivalently, $\Psi(t_1)\not\geq0$.
Let
$
\tau := \inf\{s\geq t_0 : \Psi(s) \not \geq 0 \}.
$
Then~$\tau \in  (t_0,t_1]$ and, by continuity of~$\Psi(t)$, $\Psi(\tau)\geq 0$.
Now
\begin{align}\label{eq:cco}
\Psi(\tau) &= \Psi(t_0)+\int_{t_0}^\tau \dot \Psi(s) \diff s \nonumber \\
           &= \Psi(t_0)+\int_{t_0}^\tau (rI+A(s) )\Psi(s) \diff s\nonumber \\
           &\gg 0,
	\end{align}
where we used the fact that~$\Psi(t_0)\gg 0$,
$rI+A(s)\geq0$ for almost all~$s$, and $\Psi(s)\geq0$ for all~$s \in[t_0,\tau]$.
	But this implies that there exists~$\varepsilon>0$ such~$\Psi(t) \gg  0$ 
	for all~$t\in[\tau,\tau+\varepsilon]$, and this contradicts  
	the definition of~$\tau$.
	We conclude that if~$\Theta(t_0)\gg 0$ then~$\Theta(t)\geq 0$ for all~$t\geq t_0$. 
 By continuity with respect to initial conditions,  this 
holds  also for the case~$\Theta(t_0)\geq 0$ and, in particular, 
for~$\Theta(t_0)=I$. Thus, condition~\eqref{enu:ametz}  in Lemma~\ref{lem:phiplis} implies 
condition~\eqref{enu:tnds}.

To prove the   converse implication, 
assume that~$\Theta(\tau_2,\tau_1)\geq 0$ for any pair~$(\tau_2,\tau_1)$ with~$a<\tau_1\leq \tau_2<b$.  
Fix~$t_0\in(a,b)$.
  Then for almost all~$t \in(a,b)$, 
\[
			\lim_{\varepsilon \to 0} \frac{  \Theta(t+\varepsilon ,t_0) -\Theta(t,t_0)}{\varepsilon} = A(t) \Theta(t,t_0).
\]
Multiplying on the right by~$\Theta(t_0,t)$ yields
\[
	  	\lim_{\varepsilon \to 0} \frac{  \Theta(t+\varepsilon ,t ) -I}{\varepsilon}  =A(t).
\]
Since~$\Theta(t+\varepsilon ,t ) \geq 0$, we conclude that~$A(t)$ is Metzler for almost all~$t\in(a,b)$.~\hfill{$\square$}
	
	We can now prove Thm.~\ref{thm:aprre}.
	
{\sl Proof  of Thm.~\ref{thm:aprre}.}
Suppose that~\eqref{eq:defuder}   is~TNDS.
 Then~$\Phi(t,t_0) $ is~TN for all~$a<t_0\leq t<b$.
Thus,~$\Phi^{(p)}(t,t_0)\geq 0$ for all~$a<t_0\leq t<b$ and all~$p\in\{1,\dots,n\}$, with~$\Phi^{(p)}(t_0,t_0)=I$. 
Lemma~\ref{lem:phiplis} implies that~$A^{[p]}(t)$ is Metzler for almost all~$t\in(a,b)$. In particular,~$A(t)$ and~$A^{[2]}(t)$ are Metzler for almost all~$t$
and arguing as in the second 
proof of Thm.~\ref{thm:tic} implies  that~$A(t)\in \M$ for almost all~$t$.
 
To prove the converse implication, assume that~$A(t) \in \M$ for almost all~$t$.
Pick~$p\in\{1,\dots,n\}$. Arguing as in the second 
proof of Thm.~\ref{thm:tic} implies  that~$A^{[p]}(t)$ is Metzler
 for almost all~$t$. Now
Lemma~\ref{lem:phiplis} implies that~$\Phi^{(p)}(t,t_0)\geq 0$ for all~$t\geq t_0$. 
Thus, the system is~TNDS.~\hfill{$\square$}
 
The next result provides a sufficient condition for TPDS of~\eqref{eq:defuder}.
 \begin{Theorem}\label{thm:meas}
Suppose that~$A(t) \in \M^+$  for almost all~$t\in(a,b)$ and, furthermore,
  that~$a_{ij}(t) \geq \delta>0 $ 
 for all~$|i-j|=1$ and 
  almost all~$t\in(a,b)$. Then 
the system~\eqref{eq:defuder}   is~TPDS
on~$(a,b)$.
\end{Theorem}
{\sl Proof.}
By Thm.~\ref{thm:aprre}, the system is TNDS. Pick~$1\leq p\leq n$.
To analyze~$\Phi^{(p)}$, pick~$1\leq k \leq \binom{n}{p}$, and let~$v(t)$ denote the~$k$th column of~$\Phi^{(p)}(t)$. 
Then~$\dot v(t)=A^{[p]}(t)v(t)$,
with~$v(t_0)=e^k$, where~$e^k$ is the~$k$th canonical vector in~$\R^{\binom{n}{p}}$.
 Note that every off-diagonal entry of~$A^{[p]}(t)$ is either zero or larger 
or equal to~$\delta>0$ for almost all~$t$. 
For any~$j$ the linear equation for~$\dot v_j$
implies that there exists~$c_j \in \R$ such that~$v_j(t)\geq \exp(c_j (t-t_0))v_j(t_0)$ for all~$t\geq t_0$. 
In particular, if~$v_j(\tau)>0$ at some time~$\tau$ then~$v_j(t)>0$ for all~$t\geq \tau$. Thus, $v_k(t)>0$ for all~$t\geq t_0$.
 Pick a time~$\tau \geq t_0$, and let~$s\geq 1$ denote the number of entries~$j$
such that~$v_j(\tau)>0$. Without loss of generality, assume that
 these entries are~$j=1,2,\dots, s$. Write~$\dot v=A^{[p]} v$   as 
\[
			\dot v= \begin{bmatrix} E&F\\G& H\end{bmatrix} \begin{bmatrix} v_1\\ \vdots \\
			v_s\\0    \end{bmatrix},
\]
where~$0$ denotes a vector of~$\binom{n}{p}-s$ zeros. 
Since~$A^{[p]}(t)$ is irreducible, there exists a nonzero entry in~$G$ and this entry is larger or equal to~$\delta>0$ for almost all~$t$.
 This means that at least~$s+1$ entries of~$v(t)$ are positive for all~$t>\tau$. Our assumption on~$A(t)$ implies that we can now use an inductive argument
to conclude that all the  entries of~$v(t)$ are positive for all~$t>t_0$. Since this 
holds for arbitrary~$p$ and~$k$,
we conclude that every minor of~$\Phi(t)$ is positive for all~$t>t_0$.~\hfill{$\square$}

In the special case where~$A(t)$ is continuous 
it is possible to show in a similar manner that
 necessary and sufficient condition for~TPDS on~$(a,b)$
is that~$A(t) \in \M$ for all~$t\in (a,b) $ and that none of the
 functions~$a_{i,i+1}(t),a_{i+1,i}(t)$ is zero on an interval~$[r,s]$ with
$a<r<s<b$. This result has been proved in~\cite{schwarz1970}.

The next example generalizes Example~\ref{exa:cosh}.
\begin{Example}
Let~$A(\cdot):(a,b)
\to \R^{n\times n}$ be the matrix with all entries equal to zero except for the entries on the super- and sub-diagonal that are all equal to~$t$. 
Then~\eqref{eq:defuder}  is TPDS on any interval~$(a,b)$, with~$a\geq 0$, as~$A(t)$ is continuous and the off-diagonal terms are positive except  at~$t=0$. 
  \end{Example}

\begin{figure*}[t]
 \begin{center}
  \includegraphics[scale=0.6]{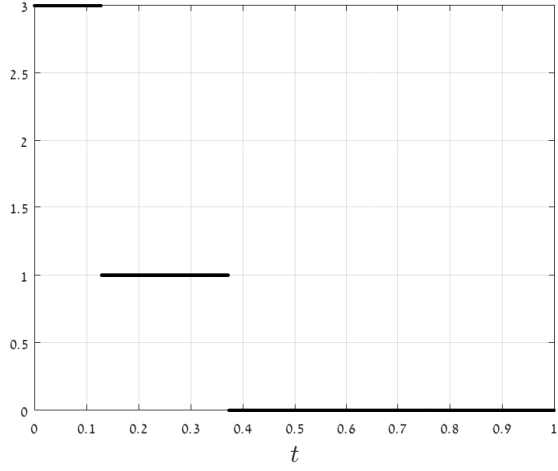}
\caption{$\sigma(z(t))$ as a function of~$t$ for the system in Example~\ref{exa:switched}. }\label{fig:imp}
\end{center}
\end{figure*}

The more general  framework that we consider  here  allows  to include
  time-varying linear systems that ``switch'' between 
 several  matrices (see, e.g.,~\cite{liberzon_book}). The next example demonstrates this.
\begin{Example}\label{exa:switched}
Consider the system
\be\label{eq:zlop}
\dot z(t)=A(t) z(t),
\ee
with 
\[
A(t):=\begin{cases}
C, & t\in[0,1/4],\\
B(t), & t\in[1/4,1/2],\\
C',&t \in[1/2,1],
\end{cases}
\]
where~$C:=\begin{bmatrix}  
-1&2& 0 &0 \\
2&-6&3&0 \\
0&5 &-1& 6\\
0&0 &4&-1  \end{bmatrix} $, $B(t)$  is the~$4\times 4 $ matrix with all entries equal to zero except for the entries 
 on the sub- and super-diagonals
 that are all equal to~$t$,  and~$C'$ is the transpose of~$C$. 
By Thm.~\ref{thm:meas}, the system is~TPDS on~$(a,b)=(0,1)$.
Fig.~\ref{fig:imp} depicts~$\sigma(z(t))$  
for the initial condition~$z(0)=\begin{bmatrix} -1& 5 &-13& 17 \end{bmatrix}'$. 
It may be seen that~$\sigma(z(t))$ is
piecewise-constant and that at any point where its  value changes it   
decreases. 
\end{Example}

%

  Thm.~\ref{thm:exp_is_tp} implies that if~$\dot z=Az$ is~TPDS on~$(a,b)$ then 
	the combined number of zeros of~$z_1(t)$ and~$z_n(t)$ on~$(a,b)$
does not exceed~$n-1$. A natural question is whether the number of zeros of other entries of the vector~$z(t)$ is  also bounded.

Consider first the case where~$A \in \M^+$ is a constant matrix. 
Then there exists~$r\geq 0$ such that~$B:=rI+A$ is~TN (by the result in  Example~\ref{exa:trid}),
irreducible, and non-singular. Thus,~$B$ is oscillatory, and applying Thm.~\ref{thm:spec} we conclude that
all the eigenvalues of~$A$
are  real   and  can be ordered as~$\lambda_1>\lambda_2>\dots >\lambda_n $. 
 Any nontrivial solution of~$\dot z=Az$ has the form
$
					z(t)=\sum_{i=1}^n c_i  \exp(\lambda_i t)   u^i  ,
$
where  the~$c_i$s are not all zero. 
This implies that \emph{any} entry of~$z(t)$ has no more than~$n-1$ isolated zeros on any time interval.

On the other-hand, the next example from~\cite{schwarz1970} shows that if~$A(t)$ is \emph{time-varying}
 then some entries of the TPDS~$\dot z=Az$ may have an unbounded number of zeros. 
\begin{Example}
Consider~\eqref{eq:linti} with~$n=3$ and~$A(t)=\begin{bmatrix} 0&1&0 \\
\frac{3}{2}-\cos(t) & 0&  \frac{3}{2}+\cos(t) \\
0&1&0
\end{bmatrix}$.
Note that~$A(t) \in \M^+$  for all~$t$, so the system is~TPDS on 
any interval~$(a,b)$.   Note that~$z(t)=\begin{bmatrix}  2+\cos(t) & -\sin(t) &-2+\cos(t)\end{bmatrix}'$ is a solution of~\eqref{eq:linti} and that~$z_2(t)$ changes sign an unbounded number of  times. 
Note also that~$z(t) \in \V$ and~$\sigma(z(t))=1$ for all~$t$.
\end{Example}

\subsection{The periodic case}
If~$\dot x=f(x)$ admits a solution~$\phi(t)$ that is periodic with period~$T$
then~$z:=\dot \phi$ satisfies~$\dot z(t)=J(\phi(t))z(t)$, and thus the matrix here is also periodic with period~$T$.
Analysis of the time-varying periodic linear system can provide considerable information on the periodic trajectory of the nonlinear system
(see, e.g.~\cite{poin_cyclic} and the references therein).

In this section, we thus consider  the TPDS~$\dot z=A(t) z$, $z(t_0)=z_0$, 
 with the additional assumption that~$A(t)$ is~$T$-periodic, i.e.
\be\label{eq:aper}
A(t)=A(t+T),\quad \text{for all } t. 
\ee

We   recall some known results from Floquet theory (see e.g.~\cite[Ch.~2]{chicone_1991}).
 For simplicity, we assume from hereon that~$t_0=0$, and   write~$\Phi(t)=\Phi(t,0)$
for the solution of~\eqref{eq:defuder} at time~$t$.
Let~$B:=\Phi(T )$. Then~$\Phi(t+T )=\Phi (t )B$ for all~$t\geq 0$. The eigenvalues 
of~$B$ are called the \emph{characteristic multipliers} of~\eqref{eq:defuder}.
If~$\alpha\in\C$, $u\in \C^n$ is an eigenvalue/eigenvector pair of~$B$ 
then~$z(t):=\Phi(t ) u $ is a solution of~$\dot z=A(t) z $, and~$z(t+T)=\alpha z(t)$. 
If~$s\in \C$ is such that~$\exp(s T)=\alpha$ then defining~$q(t):= \exp(-s t)z(t) $ yields 
$z(t)=\exp(s  t)q(t)$, and~$q(t+T)=q(t)$.

Since the system is TPDS,~$B=\Phi(T )$ is TP. Let~$\alpha_k\in \C,p^k\in \C^n$, $k=1,\dots,n$,
denote the eigenvalues and corresponding eigenvectors of~$B$.  
Thm.~\ref{thm:spec} implies that the eigenvalues   are real and satisfy
$
			\alpha_1>\alpha_2>\dots>\alpha_n>0,
$
i.e., all the characteristic multipliers are real, positive, and distinct. 
The corresponding eigenvectors 
 satisfy~$s^-(p^k)=s^+(p^k)=k-1$ for all~$k$. The next result shows that this induces   a strong structure
 on the solutions of the periodic  time-varying linear system.

\begin{Theorem}\label{thm:perioat}
Suppose that~$A(t)$  satisfies~\eqref{eq:aper}  and that~$\dot \Phi=A\Phi$ is TPDS  on~$(0,T)$. 
 Pick~$1\leq i\leq j \leq n$ and~$c_i,\dots c_j \in \R$, with~$c_i\not =0$.
 Then the  solution of~$\dot z=A(t)z$, $z(0)=\sum_{k=i}^j c_k p^k$, 
 satisfies
\be\label{eq:oiyp}
			i-1\leq \sigma(z(t))\leq  j-1,   
\ee
 for all~$t\geq 0$  except, perhaps, for up to~$j-i$ isolated points,   and
$
			  \sigma(z(t))=  i-1$   for all $t$  sufficiently large.
In particular, if~$z(0)=c_i p^i$, with~$c_i\not =0$, then
$
			\sigma(z(t))\equiv i-1$  for all $t\geq 0$.
\end{Theorem}

 Thus, the eigenvectors of~$\Phi(T)$ induce a decomposition of the state-space with respect to~$\sigma$. The monotonicity of~$\sigma$ then
restricts the possible dynamics; see~\cite{wang_trans_2015,fgwang2013}.

{\sl Proof of  Thm.~\ref{thm:perioat}.}
By Thm.~\ref{thm:spec},
\[
									i-1\leq     s^-(z(0)) \leq s^+ (z(0))       \leq j-1,
\]
and TPDS  implies that
\be\label{eq:moin}
								     s^-(z(t)) \leq s^+ (z(t))       \leq j-1.
\ee

To proceed, note that~$z(t) =\Phi(t) (\sum_{k=i}^j c_k p^k )$ implies that
\begin{align*}
				z(t+T)&=\Phi(t+T) (\sum_{k=i}^j c_k p^k )\\
						&=\Phi(t)B (\sum_{k=i}^j c_k p^k )\\
						&=\Phi(t) (\sum_{k=i}^j c_k \alpha_k p^k ),
\end{align*}
and iterating this yields
$
				z(mT)=\sum_{k=i}^j c_k (\alpha_k)^m p^k  
$,  
for any integer~$m\geq 0$.  Since~$\alpha_i> \alpha_s$ for all~$s>i$, this implies that~$\sigma( z(mT) ) =\sigma( p^i)=i-1$  for any sufficiently large~$m$.
Combining this with~\eqref{eq:moin} and using the fact that the system is~TPDS, we conclude that
\[
								   i-1\leq   s^-(z(t)) \leq s^+ (z(t))       \leq j-1, \quad \text{for all } t\geq 0, 
\]
and this proves that~\eqref{eq:oiyp} holds for all~$t\geq 0 $ 
  except, perhaps, for  up to~$j-i$ isolated points.

In the particular case, $z(0)=c_i p_i$, with~$c_i \not =0$,  \eqref{eq:oiyp} yields
\[
			i-1\leq \sigma(z(t))\leq  i-1,   
\]
for almost all~$t\geq 0$ and using monotonicity implies that~$\sigma(z(t))\equiv  i-1$.~\hfill{$\square$}

\begin{Example}
Consider the case~$n=2$ and
$
A(t)=\begin{bmatrix}
									0&1+\sin(t)\\1+\sin(t)&0 
\end{bmatrix}.
$
Note that this is~$2\pi$-periodic and yields a~TPDS on any interval~$(a,b)$. 
 The solution of~\eqref{eq:defuder} is
$
\Phi(t)=\begin{bmatrix}
									\cosh(a(t)) & \sinh(a(t)) \\  \sinh(a(t)) &  	\cosh(a(t))
\end{bmatrix},
$
where~$a(t):=1+t-\cos(t)$, so
$
				B:=\Phi(2\pi)=\begin{bmatrix}
									\cosh(2\pi) & \sinh(2\pi ) \\  \sinh(2\pi) &  	\cosh(2\pi) 
\end{bmatrix}.
$
The eigenvalues of~$B$ are~$\alpha_1=\exp( 2\pi)$, $\alpha_2=\exp(-2\pi)$ and the corresponding eigenvectors are~$ p^1=\begin{bmatrix} 1&1\end{bmatrix}'$
and~$p^2= \begin{bmatrix}-1&1\end{bmatrix}'$. 
Note that~$\sigma(p^1)= 0$ and~$\sigma(p^2)= 1$.
Consider~$z(t)$ for~$z(0)=c_1p^1+c_2p^2$ with~$c_1=1$ and~$c_2=10$. 
Then~$\sigma(z(0))=1$.
For large~$t$, $\Phi(t)\approx\frac{1}{2}\exp(a(t))J_2$, were~$J_2$ is the~$2\times 2$ matrix with all entries equal to one,   
and thus~$z(t)\approx c_1 \exp(a(t))  \begin{bmatrix} 1 & 1 \end{bmatrix}'$.
We conclude that~$\sigma(z(t))=0$ for all sufficiently large~$t$. 
\end{Example}

  Schwarz considered only  linear  time-varying systems. 
In the next section, we describe how the TPDS framework can be used to 
analyze the stability of nonlinear dynamical systems.

\section{Applications to stability analysis}\label{sec:app}

Consider  the  nonlinear  time-varying dynamical system 
\be\label{eq:gfd}
\dot  x(t)=f(t,  x(t)),
\ee
whose  trajectories   evolve  on an invariant
 set~$\Omega \subset \R^n$, that is, for any~$x_0\in\Omega$ and any~$t_0\geq 0$
a unique solution~$x(t,t_0,x_0)$ exists and satisfies~$x(t,t_0,x_0)\in\Omega$ for all~$t\geq t_0$. From here on we take~$t_0=0$.   We assume that~$\Omega$ is   compact and convex, and that~$f$ is~$C^1$ with respect to~$x$. 

\begin{Assumption}\label{ass:perj}
For any~$t\geq 0 $  and along any line~$\gamma:[0,1]\to \Omega$   
  the matrix
\be\label{eq:atpo}
A(t):=\int_0^1  J(t,\gamma(r))     \dif r
\ee
 is well-defined,  locally (essentially) bounded,  measurable, and $\dot \Phi(t) =A(t)\Phi(t)$ is~TPDS.
Here~$J(t,x):=\frac{\partial}{\partial  x} f(t,x)$ is the Jacobian of the dynamics. 
\end{Assumption}		
Note that Theorem~\ref{thm:meas} can be used to establish conditions on~$J$ guarantying the required~TPDS property.

The next result shows that under 
  Assumption~\ref{ass:perj} the system~\eqref{eq:gfd} satisfies an ``eventual monotonicity''
property.
\begin{Lemma}\label{lem:ztr}
Pick~$a,b \in \Omega$, with~$a\not = b$,  and consider the solutions~$x(t,a)$, $x(t,b)$  
  of~\eqref{eq:gfd}.  There exists a time~$s\geq 0$ such that for all~$t\geq s$
	either~$x_1(t,a)>x_1(t,b)$ or~$x_1(t,a)<x_1(t,b)$. 
	\end{Lemma}
	{\sl Proof of Lemma~\ref{lem:ztr}.}
Denote the line between the two solutions at time~$t$ by
	\[
	\gamma(  r):=rx(t,a)+(1-r)x(t,b),\quad r\in[0,1].
	\]
	Since~$\Omega$ is convex,~$\gamma ( r) \in \Omega$ for all~$t\geq 0$ and all~$r\in[0,1]$. 
Let~$z(t):=x(t,a)-x(t,b)$. Then
\begin{align*}
   \dot z(t)&=f(t,x(t,a))-f(t,x(t,b))\\
	       &=\int_0^1  \frac{d }{d r}   f(t,\gamma(r))     \dif r\\
					&=A(t)z(t),
\end{align*}
with~$A(t)$ defined  in~\eqref{eq:atpo}.
 By assumption, this is TPDS, so  Thm.~\ref{thm:exp_is_tp}
yields~$z(t) \in \V$ for all~$t$ except  for up to~$n-1$ time points. 
In particular, there exists~$s\geq 0$ such that~$z_1(t)\not = 0  $ for all~$t\geq s$.~\hfill{$\square$} 

We consider from hereon the case where~$f$ is~$T$-periodic for some~$T>0$.
\begin{Assumption}\label{ass:tper}
\[
f(t,z)=f(t+T,z), \text { for all } z\in \Omega  \text{ and } t\geq 0. 
\]
\end{Assumption} 

Note that in the particular case where~$f$ is time-invariant this property  holds for all~$T$.
 Note also  that this implies that the matrix~$A(t)$ in~\eqref{eq:atpo} is also~$T$-periodic. 
	\updt{A solution~$\gamma(t)$ of~\eqref{eq:gfd} is called a~$T$-periodic trajectory if~$\gamma(t+T)=\gamma(t)$ for all~$t$. 
}

\begin{Theorem}\label{thm:application}
If Assumptions~\ref{ass:perj} and~\ref{ass:tper}    hold then 
every solution of~\eqref{eq:gfd} converges to a~$T$-periodic trajectory. 
\end{Theorem}

This result has been derived by  \cite{periodic_tridi_smith} based on a direct analysis of the number of sign variations
in the vector of derivatives~$z(t):=\dot x(t)$.  
\updt{In the particular  case where~$f(t,x)=g(x,u)$}, with~$u(t)$  $T$-periodic, one may view~$u$ as a periodic excitation. 
Then Thm.~\ref{thm:application} implies that the system \emph{entrains} to the excitation in the sense that
every solution   converges to a periodic solution with the same period as the excitation. 
Entrainment is  important in many natural and artificial systems.
For example, proper functioning of biological organisms often requires entraining of various processes 
  to periodic excitations like the~24h
 solar day or the cell-division cycle~\citep{entrain2011,RFM_entrain}.  
Epidemics of infectious diseases often correlate with 
seasonal changes and   interventions like pulse vaccination  may also need to be periodic~\citep{epidemics_2006,2017arXiv171007321M}.

It is well-known that the linear system~$\dot x=Ax+Bu$, with~$A$ Hurwitz,  entrains to a periodic input.
More generally, contractive systems entrain (see e.g.~\cite{DBLP:journals/corr/MargaliotC17}).
However,   nonlinear systems do not necessarily entrain. 
There are examples of ``innocent looking'' nonlinear dynamical systems
that generate chaotic trajectories when excited with periodic inputs~\citep{Nikolaev145201}.

 The next example, which is a special case of a construction  from~\cite{Takac1992},
 shows that even strongly 
cooperative systems do not necessarily entrain.  
\begin{Example}
Consider the system
\be\label{eq:tak}
\dot x= \begin{bmatrix} 1&0&0&1\\ 1&1&0&0 \\ 0&1&1&0 \\ 0&0&1&1\end{bmatrix} x 
+ \begin{bmatrix}-2x_1^3+x_1 u  \\-2x_2^3-x_2 u\\-2x_3^3+x_3 u\\-2x_4^3-x_4 u
                  \end{bmatrix}.
\ee
Note that this is a strongly cooperative system, but not  a tridiagonal system because
of the feedback connection from~$x_4$ to~$x_1$.
For~$u(t)=\cos(2t)$ it is straightforward to verify that
$
				\gamma(t):= \begin{bmatrix}\cos(t)& \sin(t)& -\cos(t) &-\sin(t)\end{bmatrix}'$
is a solution of~\eqref{eq:tak}. Furthermore, it can be shown that this solution is 
locally asymptotically stable~\citep{Takac1992}. Thus,
for an excitation that is periodic with   period~$T=\pi $ there exist  
 trajectories   converging  to~$\gamma$ that are  periodic with a minimal period~$2T$. 
\end{Example}

We can now prove  Thm.~\ref{thm:application}. Pick~$a\in \Omega$.
If   the solution~$x(t,a)$  of~\eqref{eq:gfd} is $T$-periodic then there is nothing to prove. Thus,
suppose that~$x(t,a)$    is  not~$T$-periodic. 
  Then~$ x(t+T,a)$ is another solution of~\eqref{eq:gfd} that is different from~$x(t,a)$.
		Using Lemma~\ref{lem:ztr} we conclude that
	there exists 
  an integer~$m\geq 0$ such that~$x_1(k T ,a)-x_1((k+1)T,a) \not =0$ for all~$k\geq m$. Without loss of generality,
	  assume that 
\be\label{eq:xms}
					x_1(k T,a )-x_1((k+1)T,a) >0 \text{ for all } k\geq m. 
\ee
   Define the Poincar\'e
 map~$P_T:\Omega\to \Omega$ by 
$
				P_T(y):=x(T, y).
$
Then~$P_T$ is continuous, and for any integer~$k\geq 1$
 the $k$-times composition of~$P_T$ satisfies~$P_T^k(y)=x(kT,y)$. 
The  omega limit set $\omega_T:\Omega\to \Omega$ is defined by
\begin{align*}
				\omega_T(y):= & \{z\in \Omega:  \text{ there exists a sequence } n_1,n_2, \dots \\&\text{ with } n_k \to \infty 
				\text{ and }\lim_{k\to \infty} P_T^{n_k}(y)= z \} .
\end{align*}
This set is not  empty,   invariant under~$P_T$, that is,~$P_T( \omega_T(y))  =\omega_T(y)$, and
 $ P_T^n(y) \to \omega_T(y)   $ as~$n\to\infty$. 
In particular, if~$\omega_T(y)=\{q\}$ then~$P_T(q)=q$, that is, 
 the solution emanating from~$q$ is $T$-periodic. 
Thus, to prove the theorem we need to show that~$\omega_T(a)$ is a singleton. 
Assume that this is not the case. Then there exist~$p, q 
  \in \omega_T(a)$ with~$p\not = q$. 
This means  that there exist  integer sequences~$n_k\to \infty$ and~$m_k \to \infty$
such that
\[
\lim_{k\to \infty}x(n_k T ,a)=p,\quad \lim_{ k \to \infty}x(m_k T,a)=q. 
\]
Without loss of generality, we may pick $n_k < m_k < n_{k+1}$ for all $k$, which implies by~\eqref{eq:xms} that $x_1(n_k T ,a) < x_1(m_k T ,a) <x_1(n_{k+1} T ,a) $ for all $k$ sufficiently large, which passing to the limit yields~$p_1=q_1$.
In other words, any two points~$p,q\in \omega_T(a)$ have the same first coordinate. 
Consider the solutions emanating from~$p$ and from~$q$ at time zero, that is,~$x(t,p)$ and~$x(t,q)$.  
We know that  there exists an integer~$m\geq 0$ such that, say, 
\be\label{eq:grew}
						x_1(kT,p)-x_1(kT,q)>0 \text{ for all } k\geq m.
\ee
But since~$p,q\in \omega_T(a)$, $x(kT,p),x(kT,q) \in \omega_T(a)$ for all~$k$, and
  this means that~$x(kT,p)$ and~$x(kT,q)$ have the same first coordinate. This contradiction
completes the proof  of Thm.~\ref{thm:application}.~\hfill{$\square$} 


The time-invariant nonlinear   dynamical system:
\be\label{eq:gfdti}
\dot  x(t)=f( x(t))
\ee
is~$T$-periodic for all~$T>0$, so Thm.~\ref{thm:application} yields the following result.

\begin{Corollary}
Suppose that: (1)~the solutions of~\eqref{eq:gfdti}  evolve  on an invariant  compact  and convex set~$\Omega \subset \R^n$;
(2)~$f\in C^1$; and 
 (3)~the matrix~$J(x):=\frac{\partial}{\partial x}f(x) \in\M^+$   
   for all~$x\in \Omega$. Then for every~$x_0\in\Omega$ the solution~$x(t,x_0)$
  converges to an equilibrium point. 
\end{Corollary}

This is  a generalization of Smillie's theorem.
 Indeed, the proof in~\cite{smillie} is based on studying~$\sigma(z(t))$ near zeros of the~$z_i$s
and since these may be high-order zeros, Smillie used iterative differentiations and thus
had to assume  that every entry~$f_i$ of the vector field is~$(n-1)$-times differentiable~\cite[p.~530]{smillie}.


Note that the approach used by Smillie provides a sufficient condition for the
 dynamical system~$\dot z =A z$ 
 to satisfy~\eqref{eq:slop} for
all~$a<t_0<t<b$ and all~$z(t_0)\in\R^n \setminus\{0\}$, but it seems difficult to use this approach to understand if this is also a necessary condition. 
  Remark~\ref{rem:conv} above shows that  the TPDS approach solves this question.

\section{Directions for future research}\label{sec:nd}

We believe that one of the most important implications of this paper  is that 
it opens many new and interesting  research directions. We now briefly describe
several such potential directions.

 First, the elegant proofs  of Schwarz from~1970 
are based on what is now known as the theory  of cooperative  
dynamical systems. But since then this theory has been greatly developed.
Extensions  include for example   cooperative systems
in canonical form~\citep{diag_sca_smith},
the theory of monotone  (rather than cooperative) systems,  and the new notion 
of \emph{monotone control systems}~\citep{mcs_angeli_2003}. 
These extensions  can perhaps yield
new and interesting results in the context of~TPDS. 

Second, the direct proof of the~SVDP in the work of
 Smillie and others 
 is    difficult to
 generalize to other cases. Using the connection to~TPDS suggests an easier track for generalizing these results  
    to other forms of    dynamical systems, for example,  
 those with  transition \emph{operators}
that have an~SVDP.  Note that there is a large body of work on 
kernels satisfying an~SVDP   (see, e.g.~\cite{karlin_tp,Pinkus1996}).
Another possible   direction is motivated 
 by  weakening  the requirement of~TPDS to~TNDS.
Note that the proof of Thm.~\ref{thm:application} mainly uses
 the eventual monotonicity
behavior of the first entry of the solution~$z$ in a TPDS~$\dot z=A z$ 
 described in Lemma~\ref{lem:ztr}. 
In a~TNDS, this property does not hold. Yet, Schwarz
 showed that a weaker property does hold.
\begin{Lemma}\label{lem:w_ker}\citep{schwarz1970}
Suppose that~$\dot z=Az$ is a~TNDS on~$(a,b)$ and that there
 are  times~$r,s$ with~$a<r<s<b$ such that
$
z_1(r)=0$  and $z_1(s)\not =0$.
 Then
\be\label{eq:resdp}
			s^+(z(s)) \leq s^+(z(r))-1.
\ee
\end{Lemma}

Thus, TNDSs do not satisfy the eventual monotonicity described 
in Lemma~\ref{lem:ztr}, but do satisfy
 the ``non-oscillatory'' condition described in Lemma~\ref{lem:w_ker}.
The question then is whether it is possible to use this to generalize the results in the~TPDS
case to  TNDS with some additional properties.    

Another natural direction for further research is to study the 
time-discretized solutions of~TNDSs. 
 For example, consider a  
matrix~$A\in \M^+$.
	A simple  discretization of~$\dot x=Ax$ is given by 
	$
				\tilde z(k+1)=\tilde z(k)+h A \tilde z(k),
	$
with~$h>0$.
 For any~$h>0$  sufficiently small it
follows from Example~\ref{exa:trid} that the matrix~$I+hA$ is~TN. It is also nonsingular, so  
\begin{align*}
											s^-(\tilde z(k+1))  \leq  s^-(\tilde z (k) ), \;\;
											s^+(\tilde z(k+1))  \leq  s^+(\tilde z (k) ),
\end{align*}
for all~$k$. 
A similar result holds for more sophisticated discretization 
 schemes,    say when~$A(t)$ is time-varying and~$z(k+1)=z(k)+(\int A (s) \diff s) z(k)$
with   integration on an appropriate time interval.
 An interesting question is what can be deduced from these SVDPs
 on the asymptotic behavior  of the discrete-time systems.   	
	
The notion of a TPDS may be useful also for studying nonlinear time-varying, yet \emph{not} necessarily periodic, dynamical systems. Indeed, consider the system~$\dot x=f(t,x)$, and suppose that
the corresponding variational system~$\dot z=J(t,x)z$ is~TPDS. Then~$z(t) \in \V$ 
for all~$t$ except perhaps for up to~$n-1$ time points, so in particular there exists a time~$s$ such that~$z_1(t) \not =0$ for all~$t\geq s$. This means that~$x_1(t)$ is monotone for all~$t\geq s$. If the state-variables are bounded (and this is typical for example in models from systems biology) then~$x_1(t)$ converges to a limit~$e_1$.
 Similarly,~$x_n(t) \to e_n$. We may now view~$x_2(t),\dots,x_{n-1}(t)$ 
as a system of~$n-2$ state-variables with ``inputs''~$x_1(t),x_n(t)$ that converge to a constant
 value. If this system admits the converging input
converging state~(CICS) property (see, e.g.,~\cite{cics2006} and the references therein)
then we can deduce convergence to equilibrium of the entire state~$x(t)$. 

We note in passing that the fact that~$x_1(t),x_n(t)$ converge to a limit 
is interesting by itself
 especially
if these are the system outputs e.g. they feed another ``downstream'' system.

Another direction for further research is exploring the applications
of~TPDS to
  differential analysis and
contraction theory~\citep{sontag_cotraction_tutorial,forni2014,LOHMILLER1998683}. To explain this, assume that the trajectories of
\be\label{eq:pout}
\dot x=f(t,x)
\ee
 evolve on  a compact and convex  state-space~$\Omega$.
  For~$a,b\in\Omega$, let~$\gamma(r):= r a +(1-r)b$, with~$r\in [0,1]$, denote
 the line connecting~$a$ and~$b$, and let
$
					w(t,r):=\frac{\partial  }{\partial  r } x(t,t_0,\gamma(r)),
$
that is, the change in the solution at time~$t$ w.r.t. a change in the initial condition along the line~$\gamma$ at the initial time~$t_0$. 
  Then (see e.g.~\cite{entrain2011})
	\[
						\dot w(t,r)= J(t, x(t,t_0,\gamma(r))) w(t,r) .
	\]
	This is again a linear time-varying  system.  
	If it is~TPDS then one can obtain strong results on the asymptotic behaviour of~\eqref{eq:pout} using the~SVDP. This idea has already been used extensively
	by~\cite{poin_cyclic,Fusco1990} and others, but using direct analysis of the evolution of the number of sign changes in~$w$. The relation to~TPDS may lead to new results.

	Another topic for further research is based on the fact that several authors
	used a slightly different notion of
	the number of sign variations as a discrete-valued Lyapunov function.
	For a vector~$x\in\R^n$ with no zero entries 
	let~$\sigma_c(x):=|\{ i\in\{1,\dots,n \} : x_i x_{i+1}<0   \}|$, where~$x_{n+1}:=x_1$.
	This is the ``cyclic'' number of sign changes in~$x$. The results in  \cite{Fusco1990,smith_sign_changes}
	show that for some linear dynamical systems~$\sigma_c(z(t))$ can only decrease 
	along any solution~$z(t)$. The proofs are based on direct calculations. 
	Recall  that the~SVDP with respect to the ``standard'' number of sign variations~$\sigma$ 
	characterizes the  sign-regular matrices. This leads to the following question:
	when does~$A$ (and~$\exp(At)$) satisfies an SVDP
	with respect to~$\sigma_c$?

\section{Conclusions}
TN  and TP matrices enjoy a rich set of powerful
properties and have found applications in numerous fields. 
A natural question is when is the transition matrix
of   a linear dynamical system~$\dot z =A z$ TN or TP? 
This problem has been solved
 by~\cite{schwarz1970} yielding the notion of~TNDS and~TPDS.
 One important property of such systems is that for any solution~$z(t)$
the number of sign variations~$\sigma(z(t))$ is non-increasing with time. 
His approach is based on what is now known as cooperative systems theory: 
a system is~TNDS [TPDS] if all the minors of the transition matrix, that are all either zero or one at the initial time~$t_0$,  are
non-negative [positive] for all~$t>t_0$. 
However, the seminal
work of Schwarz
has been largely forgotten,  perhaps because he did not show how
 to apply these  results to analyze  \emph{non-linear} dynamical systems.

More recently, 
the number of sign changes~$\sigma(  z(t))$, where~$z:=\dot x$,
 has been used by several authors
as an integer-valued Lyapunov function for the nonlinear system~$\dot x=f(t,x)$. 
In these works, the fact that~$\sigma(  z(t))$   is non-increasing with time
has been proved by a direct and sometimes tedious analysis. 

In this paper, we reviewed these seemingly different lines of research and showed
 that the linear time-varying
system describing the evolution of~$  z$ (i.e., the variational system) 
 is
in fact~TPDS.
Our results   allow to derive important generalizations
to  several known results, while greatly simplifying the proofs.  We hope that the expository nature of this paper  
will make these fascinating topics accessible to a large audience as well as
  open  the door to many new and interesting research directions.

\noindent \emph{Acknowledgments:}
We thank the anonymous referees and the Editor 
for their constructive   comments.
The first author  is   grateful  to  
 George Weiss, Yoram Zarai, Guy Sharon,    Tsuff Ben Avraham,  Lars Gr\"une, and Thomas Kriecherbauer
for  many helpful comments. 
  
\section*{Appendix: SVDP of square sign-regular matrices}
In this Appendix, we    review the SVDP of \emph{square} sign-regular matrices.
We follow the presentation in~\cite[Ch.~V]{gk_book}, as this requires little more than
basic manipulations of determinants. 
We begin with an  auxiliary result 
that  provides information on the number of sign changes in
a vector obtained as a linear combination of~$m$ given vectors.  
\begin{Proposition}\label{prop:sv1}
Consider a set of~$m$ vectors~$u^1,\dots,u^m \in \R^n$, with~$m < n$. 
Define the matrix~$U\in \R^{n\times m }$ by
$
U:=\begin{bmatrix} u^1& u^2&\dots&u^m \end{bmatrix}.
$
The following  two conditions are equivalent:
\begin{enumerate}[(1)]
\item	\label{cond:cis} For any~$c_1,\dots,c_m\in \R$, that are not all zero,
\be\label{eq:suim}
					s^+(\sum_{k=1}^m c_i u^i)\leq m-1.
\ee
\item			\label{cond:mino}
All the minors of order~$m$ of~$U$, that is, all minors of the form 
\be\label{eq:minors}
U( i_1\;\dots\; i_m| 1\; \dots \;m), \text{ with }1\leq i_1< \dots<i_m\leq n,
\ee
 are non-zero and have the same sign. 
\end{enumerate}
\end{Proposition}

\begin{Remark}\label{rem:lininp}
Note that if~$m=n$ then condition~\eqref{cond:cis} 
always holds, whereas  condition~\eqref{cond:mino} holds iff~$U$ is non-singular.
Thus, the proposition does not hold in the case~$m=n$.  
\end{Remark}
 
		\begin{Example}
		Consider the TP matrix in Example~\ref{exa:tpandspec}. Its	
		  first two eigenvectors are
	$u^1=\begin{bmatrix} 1 & \sqrt{2} & 1\end{bmatrix}' $, and~$u^2=\begin{bmatrix} -1 & 0 & 1\end{bmatrix}' $. 
	We know from Thm.~\ref{thm:spec} that for any~$c_1,c_2\in\R$, that are not both zero,
$
			 s^+ (  c_1 u^1+c_2 u^2 )  \leq   1,
	$
	that is, condition~\eqref{cond:cis} holds. 
	In this case, $U=\begin{bmatrix}  1 & -1\\\sqrt{2} & 0\\1&1   \end{bmatrix}$, and thus
	the minors in~\eqref{eq:minors} are
	$
	\det(\begin{bmatrix}  1 & -1\\\sqrt{2} & 0 \end{bmatrix})$,
	$\det(\begin{bmatrix}  1 & -1\\1 & 1 \end{bmatrix})$,
  $\det( \begin{bmatrix}  \sqrt{2} & 0\\1 & 1 \end{bmatrix})$.
	These  are all positive, so  condition~\eqref{cond:mino} also holds.
		\end{Example}

{\sl Proof of Prop.~\ref{prop:sv1}.}
It is enough to prove the result  for the case~$n=m+1$ (the proof when~$n>m+1$ is similar). 

We first show that
condition~\eqref{cond:mino} implies condition~\eqref{cond:cis}. 
Suppose that all the minors in~\eqref{eq:minors}  are non-zero and have the same  sign.
Pick~$c_1,\dots,c_m\in \R$, that are not all zero,  and let~$u:=\sum_{i=1}^m c_i u^i \in \R^{m+1}$. 
We need to show that~$s^+(u)\leq m-1$. 
Seeking a contradiction, assume that~$s^+(u)\geq m $. This means that
\be\label{eq:uii}
u_iu_{i+1}\leq 0  \text{ for  }i=1,\dots,m .
\ee
Furthermore, condition~\eqref{cond:mino} implies that at least one of the first~$m$  
entries of~$u$ is not zero. Consider the square
 matrix~$V:=\begin{bmatrix}u&u^1&\dots & u^m \end{bmatrix}$. 
Expanding its determinant along the first column and using~\eqref{eq:uii} and condition~\eqref{cond:mino}
implies that~$\det(V) \not =0$. But, the first column of~$V$ is a linear combination of the other columns, so~$\det(V)=0$. This contradiction completes the proof that
condition~\eqref{cond:mino} implies condition~\eqref{cond:cis}. 

To prove the converse implication, assume that
condition~\eqref{cond:cis} holds.
This implies in particular that for any~$c_1,\dots,c_m$, that are not all zero,
$\sum_{k=1}^m c_i u^i\not =0$, 
so the~$u^i$s are linearly independent.  
 Assume that one of the 
minors in~\eqref{eq:minors} is zero. Then there exist~$c_1,\dots,c_m$, not all zero,
 such that~$\sum_{i=1}^m c_k u^k$ has~$m$ zero entries. Since~$n=m+1$ this means that~$s^+(\sum_{i=1}^m c_k u^k)\geq m$ and this is a contradiction.
 We conclude that all the minors of order~$m$ are not zero. 
To prove that they all have the same sign it is enough to prove that 
all the minors  
\[
				d_k:=U(  1\;\dots\; k-1\; k+1 \;\dots\; m+1  | 1\;\dots \;m)  , 
\]
with~$ k=1,\dots,m+1,$ have the same sign.
Fix~$1\leq i<j \leq m+1$. Define a vector~$z \in\R^{m+1}$ by
\[
z_k:=\begin{cases} (-1)^{i-1}d_j & \text{if } k=i, \\
                   (-1) ^j d_i & \text{if }  k=j,\\
									0 & \text{otherwise.} \end{cases}
\] 
Consider the square
 matrix~$V:=\begin{bmatrix} u^1&\dots & u^m &z \end{bmatrix}$. 
Expanding its determinant along the last column 
yields~$\det(V)   =0$.
Thus, there exist~$c_1,\dots,c_{m+1}$, not all zero, such that
\begin{align}\label{eq:drfgb}
v&:=\sum_{k=1}^{m}c_k u^k \\&=-c_{m+1} z .
\nonumber 
\end{align}
  If~$c_1=\dots=c_m=0$ then this gives~$z=0$, but 
this is a contradiction as~$d_j,d_i\not = 0 $. Thus, at least one of~$c_1,\dots,c_m$ is not zero. If~$c_{m+1}=0$ then~\eqref{eq:drfgb} yields~$\sum_{k=1}^{m}c_k u^k=0$ and this contradicts condition~\eqref{cond:cis}. 
We conclude that~$c_{m+1}\not =0$, and we assume from here on that~$c_{m+1}=1$, so
\be\label{eq:psftpp}
v=-z.   
\ee
Let~$\alpha:=\sgn(d_j)\in\{-1,1\}$ and 
let~$\bar v$ be the vector~$v$ but with every zero entry~$v_p$
 replaced by~$(-1)^p  \alpha$. Then  by the definition of~$s^+$,
$s^+(v)\geq s^+(\bar v)$.
Using~\eqref{eq:psftpp} implies that the entries of~$\frac{1}{\alpha}\bar v$ are
\begin{align*}
 (-1)^1   & ,\; (-1)^2,   \; \dots \;
(-1)^{i-1}  , \; \;\frac{ (-1)^i}{\alpha}  d_j ,\; \; (-1)^{i+1}  ,   \dots ,\\&
(-1)^{j-1}   ,\;\;   \frac{(-1)^{j-1}}{\alpha} d_i ,\; \; (-1)^{j+1}   ,    \dots  , (-1)^{m+1}   .
\end{align*}
If~$d_i,d_j$ have different signs then we see that~$s^+(\bar v)\geq m$,
so~$s^+(v)\geq m$. 
This contradiction completes the proof that
condition~\eqref{cond:cis} implies
condition~\eqref{cond:mino}.~\hfill{$\square$}

We can now state the main result in this 
Appendix.

\begin{Theorem}\label{thm:sv2sq}
Consider a set of~$n$ \emph{linearly  
independent} vectors~$a^1,\dots,a^n \in \R^n$. 
The following  two conditions are equivalent:
\begin{enumerate}[(1)]
\item	\label{cond:ccombsq} For any vector~$c\in\R^n\setminus\{0\}$,
\be\label{eq:spsmmsq}
					s^+(\sum_{k=1}^n c_i a^i)\leq s^-(c).
\ee
\item			\label{cond:minoofasq}
The    matrix
$
A:=\begin{bmatrix} a^1& a^2&\dots&a^n \end{bmatrix} 
$
is SSR. 
\end{enumerate}
\end{Theorem}
 
 Clearly, Thm.~\ref{thm:sv2sq} is equivalent to Thm.~\ref{thm:regu}.

{\sl Proof of Thm.~\ref{thm:sv2sq}.}
Assume that condition~\eqref{cond:ccombsq}  holds. 
Pick 
 $1\leq p\leq n$.
We will show that all minors  of  order~$p$ of~$A$  are non-zero and have the same sign.
If~$p=n$ then this holds because the~$a^i$s are linearly independent. 
Thus,  consider   the case~$p<n$. 
Pick~$p+1$ indices~$1\leq k_1<k_2<\dots<k_p<k_{p+1}\leq n$. 
For any~$c\in\R^n$, define~$\bar c\in \R^n$
by
\[
				\bar c_i:=\begin{cases} c_i & \text{if }  i\in\{k_1,\dots,k_p\}, \\
				                          0 & \text{otherwise}. 
																	\end{cases} 
\] 
Then~$s^-(\bar c)\leq p-1$, as~$\bar c$ has no more than~$p$ non-zero entries. 
Let~$\bar a:=\sum_{i=1}^n \bar{c}_i a^i =\sum_{j=1}^p c_{k_j}a^{k_j}$. 
Applying condition~\eqref{cond:ccombsq} 
 implies that for any~$\bar c \not = 0$, $s^+(\bar a)\leq p-1$. 
This means that the set~$\{a^{k_1},\dots,a^{k_p}\}$ 
satisfies condition~\eqref{cond:cis}  in  
Prop.~\ref{prop:sv1}. Thus, all minors of the form
\be\label{eq:aaminors}
A( i_1\;\dots\; i_p| k_1\; \dots \;k_p), \text{ with }1\leq i_1 <\dots<i_p\leq n,
\ee
 are non-zero and have the same sign. Denote this sign by~$\varepsilon(k_1, \dots ,k_p)$. It remains to show that this sign
depends  on~$p$, but not on the particular choice of~$k_1,\dots,k_p$.
Pick~$v\in\{1,\dots,p\}$.
We will show that
\begin{align}\label{eq:vere}
										\varepsilon&(k_1,\dots,k_{v-1},k_{v+1},\dots,k_{p+1}) \nonumber \\&= 
										\varepsilon(k_1,\dots,k_{v},k_{v+2},\dots,k_{p+1 }).
\end{align}
To do this, define~$p$ vectors~$  \bar a^{k_1},\dots,\bar a^{k_v},\bar a^{k_{v+2}}
,\dots,\bar a^{k_{p+1}}   $ by
$\bar a^{k_i}:=a^{k_i}$ for~$i\in\{ 1,\dots,v-1,v+2,\dots,p+1\}$, 
and~$\bar a^{k_v} :=d_v a^{k_v}+d_{v+1} a^{k_{v+1}}$, where~$d_v,d_{v+1}>0$.
Pick~$\bar c_1,\dots,\bar c_v,\bar c_{v+2},\dots,\bar c_{p+1}$,
 that are not all zero, and 
let
\be\label{eq:deftra}
a:=\sum_ {   \mycom{i=1}{i\not = v+1} }^{p+1} \bar  c_i \bar a^{k_i}. 
\ee
Then 
\be\label{eq:aprt}
				a=\sum_{i=1 }^{p+1} g_i  a^{k_i},
\ee
with~$g_v:=\bar c_v d_v$, $g_{v+1}:=\bar  c_{v} d_{v+1}$, and~$g_i=\bar c_i$ for
 all other~$i$. Let~$g:=\begin{bmatrix} g_1&\dots&g_{p+1}\end{bmatrix}'$.
Note that~$g\not =0$, and that since~$g_v g_{v+1}=(\bar c_v)^2 d_vd_{v+1}\geq 0$, 
$s^-(g)\leq p-1$. 
Applying condition~\eqref{cond:ccombsq}  to~\eqref{eq:aprt}
yields~$s^+(a)\leq s^-(g)=p-1$, that
 is,
\[
s^+(\sum_ {\mycom{i=1} {i\not = v+1} }^{p+1} \bar  c_i \bar a^{k_i})\leq p-1.
\]
 Let~$\bar A \in\R^{n\times p}$ be the matrix
\begin{align*} 
\bar A:&=\begin{bmatrix} \bar a^{k_1}&\dots&\bar a^{k_{v-1}}& \bar a^{k_v}&
\bar a^{k_{v+2}}&\dots&\bar a^{k_{p+1}} \end{bmatrix}
\\&= \begin{bmatrix}   a^{k_1}&\dots&  a^{k_{v-1}}& d_v a^{k_v}+d_{v+1}a^{k_{v+1}}&
  a^{k_{v+2}}&\dots&  a^{k_{p+1}} \end{bmatrix}.
\end{align*}
Applying Prop.~\ref{prop:sv1} to the set of~$p$ 
vectors~$\bar a^{k_1}$, $\dots$, $\bar a^{k_v}$, $\bar a^{k_{v+2}}$, $\dots$, $\bar a^{k_{p+1}}\in\R^n$ we conclude that
 all minors
\begin{align}\label{eq:twomin}
&\bar  A  (    i_1,\dots,i_p  |  k_1,\dots,k_v,k_{v+2},\dots,k_{p+1} )\nonumber \\&=
d_v A(  i_1,\dots,i_p  |  k_1,\dots,k_v,k_{v+2},\dots,k_{p+1} )\nonumber \\&
+d_{v+1} A(  i_1,\dots,i_p  |  k_1,\dots,k_{v-1},k_{v+1},\dots,k_{p+1} )
\end{align}
are non-zero. This holds for  all~$d_v,d_{v+1}>0$,
so  the two minors on the right hand side of~\eqref{eq:twomin} have the same sign. 
Since this is true for all~$v\in\{1,\dots,p\}$,
we conclude that the sign~$
	\varepsilon(k_1,\dots,k_p) 
$
does not change if we change any one of the indices~$k_i$, and thus it is independent of
the choice of~$k_1,\dots,k_p$. 
This completes the proof that condition~\eqref{cond:ccombsq} implies condition~\eqref{cond:minoofasq}.

To prove the converse implication, assume that  
$
A =\begin{bmatrix} a^1&  \dots&a^n \end{bmatrix} 
$
is~SSR.  Pick~$c\in \R^n\setminus\{0\}$, and let~$p:=s^-(c)$.
If~$p=n-1$ then clearly~$s^+(\sum_{i=1}^n c_i a_i) \leq  s^-(c)$. Consider the case~$p <  n-1$. 
We may assume that the first non-zero entry of~$c$ is positive. 
Then~$c$ can be decomposed  into~$p+1$ groups:
$
(c_1,\dots,c_{v_1})$, $(c_{v_1+1},c_{v_1+2},\dots,c_{v_2})$, $\dots$,
$(c_{v_p+1},c_{v_p+2},\dots,c_{v_{p+1}})$,
where
 $c_1,\dots,c_{v_1}\geq 0$ (with at least one of these entries positive);
$c_{v_1+1}<0$, $c_{v_1+2},\dots, c_{v_2}\leq 0$, $c_{v_2+1}>0$, and so on. 
Define vectors~$u^1,\dots, u^{p+1} \in \R^{n}$ by
\[
				u^1:=\sum_{k=1}^{v_1}|c_k|a^k, \; u^2:= \sum_{k=v_1+1}^{v_2}|c_k|a^k,\dots.
\]
Then
\be\label{eq:asqre}
a:=\sum_{k=1}^{n}c_k a^k = u^1-u^2+u^3-\dots+(-1)^p u^{p+1}.
\ee
Note that every~$u^i$ is a non-negative and non-trivial sum of a consecutive set of~$a^k$s. 
Let
$
U:=\begin{bmatrix} u^1 &  \dots & u^{p+1}\end{bmatrix} \in\R^{n \times(p+1)}. 
$ Note that~$n>p+1$. 
The~SSR of~$A$ implies that 
all the minors of order~$(p+1)$ of~$U$ are non-zero and have the same sign.  
 Applying Prop.~\ref{prop:sv1} to~\eqref{eq:asqre}
yields~$s^+(a) \leq p $, so~$s^+(a) \leq  s^-(c)$.~\hfill{$\square$}

 \small

{\bf Michael Margaliot} received the BSc (cum laude) and MSc degrees in
 Elec. Eng. from the Technion-Israel Institute of Technology-in
 1992 and 1995, respectively, and the PhD degree (summa cum laude) from Tel
 Aviv University in 1999. He was a post-doctoral fellow in the Dept. of
 Theoretical Math. at the Weizmann Institute of Science. In 2000, he
 joined the Dept. of Elec. Eng.-Systems, Tel Aviv University,
 where he is currently a Professor and Chair. His  research
 interests include the stability analysis of differential inclusions and
 switched systems, optimal control theory, fuzzy control, computation with
 words, Boolean control networks, contraction theory, and systems biology.
 He is co-author of \emph{New Approaches to Fuzzy Modeling and Control: Design and
 Analysis}, World Scientific,~2000 and of \emph{Knowledge-Based Neurocomputing}, Springer,~2009. 
 He  served
as  an Associate Editor of~\emph{IEEE Transactions on Automatic Control} during 2015-2017.

{\bf Eduardo Sontag} received his Licenciado degree from the Mathematics Department at the University of Buenos Aires in 1972, and his Ph.D. (Mathematics) under Rudolf E. Kalman at the University of Florida, in 1977. From 1977 to 2017, he was with the Department of Mathematics at Rutgers, The State University of New Jersey, where he was a Distinguished Professor of Mathematics as well as a Member of the Graduate Faculty of the Department of Computer Science and the Graduate Faculty of the Department of Electrical and Computer Engineering, and a Member of the Rutgers Cancer Institute of NJ. In addition, Dr. Sontag served as the head of the undergraduate Biomathematics Interdisciplinary Major, Director of the Center for Quantitative Biology, and Director of Graduate Studies of the Institute for Quantitative Biomedicine. In January 2018, Dr. Sontag was appointed as a University Distinguished Professor in the Department of Electrical and Computer Engineering and the Department of
  BioEngineering at
  Northeastern University.  Since 2006, he has been a Research Affiliate at the Laboratory for Information and Decision Systems, MIT, and since 2018 he has been a member of the Faculty in the Program in Therapeutic Science at Harvard Medical School.

His major current research interests lie in several areas of control and dynamical systems theory, systems molecular biology, cancer and immunology, and computational biology.
He has authored over five hundred research papers and monographs and book chapters in the above areas with over 43,000 citations and an h-index of 89.   He is in the Editorial Board of several journals, including: IET Proceedings Systems Biology, Synthetic and Systems Biology, International Journal of Biological Sciences, and Journal of Computer and Systems Sciences, and is a former Board member of SIAM Review, IEEE Transactions in Automatic Control, Systems and Control Letters, Dynamics and Control, Neurocomputing, Neural Networks, Neural Computing Surveys, Control-Theory and Advanced Technology, Nonlinear Analysis: Hybrid Systems, and Control, Optimization and the Calculus of Variations. In addition, he is a co-founder and co-Managing Editor of the Springer journal MCSS (Mathematics of Control, Signals, and Systems).

He is a Fellow of various professional societies: IEEE, AMS, SIAM, and IFAC, and is also a member of SMB and BMES. He has been Program Director and Vice-Chair of the Activity Group in Control and Systems Theory of SIAM, and member of several committees at SIAM and the AMS, including Chair of the Committee on Human Rights of Mathematicians of the latter in 1981-1982. He was awarded the Reid Prize in Mathematics in 2001, the 2002 Hendrik W. Bode Lecture Prize and the 2011 Control Systems Field Award from the IEEE, the 2002 Board of Trustees Award for Excellence in Research from Rutgers, and the 2005 Teacher/Scholar Award from Rutgers.

\bibliographystyle{abbrvnat}
\bibliography{tpds_bib}
 \end{document}